\documentclass[journal ]{new-aiaa}
\usepackage[utf8]{inputenc}
\usepackage{textcomp}

\usepackage{graphicx}
\usepackage{subfigure}
\usepackage{amsmath}
\usepackage[version=4]{mhchem}
\usepackage{siunitx}
\usepackage{longtable,tabularx}
\usepackage{makecell}
\usepackage{hyperref}
\title{Analysis of pitchfork bifurcations and symmetry breaking in the elliptic restricted three-body problem}

\begin{document}
\author[1,2]{Haozhe Shu}
\affil[1]{Mathematical Institute, Tohoku University, Sendai, 980-8578, Japan}
\affil[2]{Advanced Institute for Material Research, Tohoku University, Sendai, 980-8577, Japan}
\author[2]{Mingpei Lin\footnote{Corresponding author: lin.mingpei.d2@tohoku.ac.jp}}
\maketitle

\begin{abstract}
A unified framework is proposed to quantitatively characterize pitchfork bifurcations and associated symmetry breaking in the elliptic restricted three-body problem (ERTBP). 
It is known that planar/vertical Lyapunov orbits and Lissajous orbits near the collinear libration points undergo pitchfork bifurcations with varying orbital energy. These bifurcations induce symmetry breaking, generating bifurcated families including halo/quasi-halo orbits, axial/quasi-axial orbits, and their corresponding invariant manifolds. 
Traditional semi-analytical methods for constructing halo orbits, based on resonant bifurcation mechanisms, have obstacles in fully exploiting the intrinsic symmetry breaking characteristics in pitchfork bifurcations. 
In this paper, a unified trigonometric series-based framework is proposed to analyze these bifurcated families from the perspective of coupling-induced bifurcation mechanisms.
By introducing a coupling coefficient and various bifurcation equations into the ERTBP, different symmetry breaking is achieved when the coupling coefficient is non-zero.
This unified semi-analytical framework captures bifurcations of both periodic/quasi-periodic and transit/non-transit orbits. Furthermore, it reveals that pitchfork bifurcation solutions in the ERTBP fundamentally depend solely on the orbital eccentricity and three amplitude parameters of the system's degrees of freedom, governing both the elliptic direction and the hyperbolic one.
\end{abstract}

\section{Introduction}
The spatial circular restricted three-body problem (CRTBP), which describes the motion of an infinitesimal body (spacecraft or asteroid) under the gravitational influence of two primaries in circular orbits, serves as a foundational model in celestial mechanics \cite{GomezMondelo,GomezMasdemontMondelo}.
As a first natural extension, the elliptic restricted three-body problem (ERTBP) generalizes this framework by allowing the primaries to evolve along small-eccentricity Keplerian orbits. The ERTBP introduces explicit time dependence and results in the absence of a first integral \cite{Broucke,Ovenden}.
This modification transforms the five stationary Lagrangian points in CRTBP into "instantaneous equilibrium positions" with periodic oscillations.
Compared to the corresponding CRTBP, the ERTBP provides a more precise and essential framework for modeling real-world celestial systems, particularly in scenarios involving planetary resonances, asteroid dynamics, and planet-moon systems \cite{parker2014low}. Deep insights into these dynamical structures also enable optimized spacecraft trajectories, including libration point orbit transfers \cite{Peng_2,Jorba_Cusco}, station-keeping strategies \cite{Shirobokov,Gurfil}, and resonant flyby for deep-space exploration \cite{Antoniadou}.

The dynamics near the collinear libration points in the restricted three-body problem (RTBP) have been extensively studied through both numerical and analytical approaches.
For numerical perspectives, Peng et al. \cite{Peng_1} generated multi-revolution halo orbits via continuation methods and multi-segments optimization methods. Initializing with periodic orbits in the CRTBP, Ferrari et al. \cite{Ferrari} succeeded in finding periodic orbits in the ERTBP through a differential correction algorithm. Paez et al. \cite{Paez_Guzzo} classified transit orbits in the ERTBP through the Floquet-Birkhoff normalization approach. Recently, Jorba et al. \cite{Jorba_chaos} systematically analyzed Hilda asteroids, which are related to a 3:2 orbital resonance with Jupiter and provided a comparison work between the CRTBP and the ERTBP. 

For analytical perspectives, high-precision orbit approximations are essential for characterizing local dynamics, offering both physical insights and initial guesses for numerical methods.
To derive approximate high-order analytical solutions, two prominent techniques have been developed: the Lindstedt-Poincaré perturbation method and the Hamiltonian normal form method.
Farquhar \cite{Farquhar} pioneered the concept of halo orbits using a frequency control scheme. Building on this foundation, 1:1 resonant bifurcation mechanisms were developed to systematically obtain halo families semi-analytically. Based on the synchronization construction of the in-plane and out-of-plane oscillators, Richardson \cite{Richardson} constructed a third-order analytical solution for halo orbits by incorporating a correction term to adjust the out-of-plane frequency. Further advancement was made by Masdemont \cite{Masdemont}, who derived high-order series solutions for invariant manifolds in the CRTBP. Extending these results to the elliptic case, Lei et al. \cite{Lei} completed high-order semi-analytical computations for the ERTBP.

In contrast to the Lindstedt-Poincaré method, the Hamiltonian normal form offers a distinct yet powerful approach for analyzing the local dynamics. 
Within this framework, Jorba and Masdemont \cite{Jarba_Masdemont} obtained a qualitative description of the local phase space by analyzing the reduced Hamiltonian in the CRTBP. Paez and Guzzo \cite{Paez_Guzzo_Ana} presented a semi-analytical construction of halo orbits and halo tubes in the elliptic model using the Floquet-Birkhoff resonant normal form. In the most recent, based on the Lie transform, Celletti et al. \cite{Celleti} presented an explicit resonant normal form, enabling analytical investigations of planar/vertical Lyapunov orbits and halo orbits in the ERTBP.

Even though, based on the resonant framework, the first-level bifurcations associated with periodic orbits have been thoroughly analyzed by using semi-analytical approaches, the indirect relationship between 1:1 resonant bifurcation mechanisms and potential bifurcations associated with quasi-periodic orbits and transit/non-transit orbits motivates some other insights. Lin et al. \cite{Lin_1,Lin_2} initially addressed coupling-induced mechanisms in the CRTBP. Without relying on the classical constructions of resonant modification, a so-called coupling modification was comparably introduced to finally obtain a comprehensive semi-analytical construction of the local phase space in the CRTBP.  

In this paper, we present a systematic analysis of pitchfork bifurcations and associated symmetry breaking near the collinear libration points in the ERTBP, within a unified trigonometric series-based framework. 
Known as a three degree-of-freedom (DOF) non-autonomous Hamiltonian system, the governing differential equations of the ERTBP totally exhibit three distinct types of symmetries.
Pitchfork bifurcations arise near the collinear libration points, accompanied by the breaking of certain characteristic symmetries in the non-bifurcated solutions \cite{Doedel}.
To achieve the corresponding symmetry breaking from non-bifurcated solutions, we introduce bifurcation equations and coupling coefficients along specific directions in the dynamical model of the ERTBP.  
These coefficients are derived by solving the coupled bifurcation equations, while non-zero coupling coefficients lead to symmetry breaking in the original solutions, triggering pitchfork bifurcations.

Unlike previous studies of semi-analytical construction around the collinear libration points, the hyperbolic part of the solution is formulated in a trigonometric form defined on the complex plane. Taking advantage of the symmetries of the ERTBP, this unified trigonometric framework aids in simplifying the specific perturbation computation. Moreover, by introducing a complex-valued amplitude parameter in describing the motion in the hyperbolic direction, the orbital bifurcation appeared in the hyperbolic part of the solution becomes able to be quantified.
Within this framework, despite the bifurcation analysis associated with center manifolds, the analysis of the constructed bifurcation equation also provides deep insights into the bifurcation of the transit/non-transit orbits. All types of bifurcated solutions in the form of trigonometric series are shown to depend solely on the eccentricity and three amplitudes corresponding to the system's DOFs, where the critical conditions are derived explicitly.

The remainder of the paper is organized as follows. In Section \ref{sec2}, the dynamical model of the ERTBP is introduced in a pulsating-rotating frame. In Section \ref{sec3}, we present a unified trigonometric series-based semi-analytical construction for the bifurcated orbits. A quantitative analysis of pitchfork bifurcations and symmetry breaking is proposed in Section \ref{sec:results}. Finally, Section \ref{sec5} offers conclusions.
\section{Dynamical model}\label{sec2}
In this section, the dynamical model of the ERTBP is introduced. The model describes the motion of an infinitesimal particle in the gravitational field of two primaries. By neglecting the attraction influence of the particle on the primaries, the motion of two primaries can be described by the Kepler orbits around their common centroid. The classical pulsating-rotating frame is employed to simplify the formulation of the dynamical model. Specifically, positioning the origin at the centroid of the two primaries, the orientation of the $X$-axis is given by the line that goes from the smaller primary to the larger primary, while the $Z$-axis has the orientation determined by the angular motion of the primaries. $Y$-axis completes the right-handed coordinate system. In this case, the normalized coordinates for the smaller and the larger primary are $(\mu-1,0,0)$ and $(\mu,0,0)$ respectively, where $\mu=m_2/(m_1+m_2)$ is the system parameter, representing the mass ratio of the smaller celestial body to the sum of the masses of the two bodies. The governing differential equations describing the motion of the infinitesimal particle in the normalized pulsating-synodic frame are as follows \cite{Szebehely}:
\begin{equation}\label{ERTBP1}
    \begin{aligned}
        X^{\prime\prime}-2Y^{\prime}=\frac{\partial \Omega}{\partial X},\\
        Y^{\prime\prime}+2X^\prime=\frac{\partial \Omega}{\partial Y},\\
        Z^{\prime\prime}+Z=\frac{\partial\Omega}{\partial Z},
    \end{aligned}
\end{equation}
where
$$\Omega(X,Y,Z,f)=\frac{1}{1+e\cos{f}}[\frac{1}{2}(X^2+Y^2+Z^2)+\frac{1-\mu}{r_1}+\frac{\mu}{r_2}+\frac{1}{2}\mu(1-\mu)]$$
represents the potential in the ERTBP. Here, $e$ is orbital eccentricity and $f$ represents the true anomaly of the secondary on the elliptic orbit. The derivatives of coordinates are defined as
$
X^\prime:=\frac{dX}{df},\ Y^{\prime}:=\frac{dY}{df},\ Z^\prime:=\frac{dZ}{df}.
$
$r_1$ and $r_2$ denote the instantaneous distances from the spacecraft to the massive and secondary primaries satisfying
\begin{equation}
    \begin{aligned}
        {r_1}^2&=(X-\mu)^2+Y^2+Z^2,\\
        {r_2}^2&=(X+1-\mu)^2+Y^2+Z^2.
    \end{aligned}
\end{equation}
The non-autonomous nature of governing equations of the ERTBP introduces complicated time-dependent perturbations. When the eccentricity $e=0$, this dynamical model reduces to the well-known autonomous CRTBP. Similar to the CRTBP, the ERTBP (\ref{ERTBP1}) possesses two independent kinds of symmetries given by \cite{Lei}:
\begin{equation}\label{sy}
    \begin{aligned}
        S_1:(f,X,Y,Z,X^\prime,Y^\prime,Z^\prime)\longleftrightarrow (f,X,Y,-Z,X^\prime,Y^\prime,-Z^\prime),\\
        S_2:(f,X,Y,Z,X^\prime,Y^\prime,Z^\prime)\longleftrightarrow (-f,X,-Y,Z,-X^\prime,Y^\prime,-Z^\prime),
    \end{aligned}
\end{equation}
where a third symmetry of the ERTBP can be obtained by directly composing the $S_1$-type and $S_2$-type symmetries:
$$
S_3:(f,X,Y,Z,X^\prime,Y^\prime,Z^\prime)\longleftrightarrow(-f,X,-Y,-Z,-X^\prime,Y^\prime,Z^\prime).
$$
For instance, suppose that a curve $(X(f),Y(f),Z(f))$ solves (\ref{ERTBP1}). It is clear that its reflection about the $(X,Y)$ plane, $(X(f),Y(f),-Z(f))$ is also a solution to (\ref{ERTBP1}). Moreover, after time reversal, its reflection about the $(X,Z)$ plane, $(X(-f),-Y(-f),Z(-f))$, also satisfies the governing differential equations. 
In (\ref{ERTBP1}), there are five Lagrange points pulsating in the synodic coordinate system. Three of these are collinear libration points, while the remaining two are triangular libration points. Inheriting the notation in the CRTBP, we denote the collinear libration points by $L_i\ (i=1,2,3)$. By adopting the following transformations of coordinates reference, the origin of (1) can be relocated to the collinear libration points:
\begin{equation}
\begin{aligned}
    &X=-\gamma_ix+\mu+a_i,\ Y=-\gamma_iy,\ Z=\gamma_iz,\ i=1,2;\\
    &X=\gamma_ix+\mu+\gamma_i,\ Y=\gamma_iy,\ Z=\gamma_iz,\ i=3,
\end{aligned}
\end{equation}
where $a_1=-1+\gamma_1$ and $a_2=-1-\gamma_2$. $\gamma_i$ denotes the instantaneous distance between the libration point $L_i$ and its closest primary. It is known that the value of $\gamma_i$ is determined by the unique positive root of Euler's quintic equation \cite{Celletti_book,Jarba_Masdemont}:
\begin{equation}
\begin{aligned}
    &{\gamma_i}^5\mp(3-\mu){\gamma_i}^4+(3-2\mu){\gamma_i}^3-\mu{\gamma_i}^2\pm2\mu\gamma_i-\mu=0,\ i=1,2;\\
    &{\gamma_i}^5+(2+\mu){\gamma_i}^4+(1+2\mu){\gamma_i}^3-(1-\mu){\gamma_i}^2-2(1-\mu)\gamma_i-(1-\mu)=0,\ i=3,
\end{aligned}
\end{equation}
where the upper and lower signs correspond to $L_1$ and $L_2$, respectively. In the transformed coordinate system, the governing equations (\ref{ERTBP1}) become
\begin{equation}
    \label{ERTBP2}
    \begin{aligned}
        x^{\prime\prime}-2y^\prime=\frac{1}{{\gamma_i}^2}\frac{\partial \Omega}{\partial x},\\
        y^{\prime\prime}+2x^\prime=\frac{1}{{\gamma_i}^2}\frac{\partial \Omega}{\partial y},\\
        z^{\prime\prime}=\frac{1}{{\gamma_i}^2}\frac{\partial \Omega}{\partial z},
    \end{aligned}
\end{equation}
with $$\Omega(x,y,z,f)=\frac{1}{1+e\cos{f}}[\frac{1}{2}((\mu-1- \gamma_i (x\mp1))^2+{\gamma_i}^2y^2+{\gamma_i}^2z^2)+\frac{1-\mu}{r_1}+\frac{\mu}{r_2}+\frac{1}{2}\mu(1-\mu)],\ i=1,2.$$ 
The upper sign is for $L_1$ and lower one for $L_2$. To obtain a high-order semi-analytical construction for bifurcated orbits near the collinear libration points, the right-hand side of (\ref{ERTBP2}) is expanded into a recurrent form \cite{Lei}:
\begin{equation}
    \label{ERTBP3}
    \begin{aligned}
    x^{\prime\prime}-2y^\prime-(1+2c_2)x=&\sum\limits_{i\geq1}[(1+2c_2)x(-e)^i\cos^if]\\
    &+\sum\limits_{i\geq0}\{(-e)^i\cos^if[\sum\limits_{n\geq2}c_{n+1}(n+1)T_n(x,y,z)]\},\\
    \ \ y^{\prime\prime}+2x^\prime+(c_2-1)y=&\sum\limits_{i\geq1}[(1-c_2)y(-e)^i\cos^if]\\
    &+\sum\limits_{i\geq0}\{(-e)^i\cos^if[y\sum\limits_{n\geq2}c_{n+1}R_{n-1}(x,y,z)]\},\\
    \ \ \ \ \ \ \ \ \ \ \ \ \ \ \ \ \ \ \ \ \ \ z^{\prime\prime}+c_2z=&\sum\limits_{i\geq1}[-c_2 z(-e)^i\cos^if]\\
    &+\sum\limits_{i\geq0}\{(-e)^i\cos^if[z\sum\limits_{n\geq2}c_{n+1}R_{n-1}(x,y,z)]\}.
    \end{aligned}
\end{equation}
Here, $\{T_n\}$ and $\{R_n\}$ are sequences of homogeneous polynomials defined by the recurrence relations 
\begin{equation}
    \begin{aligned}
        T_n=\frac{2n-1}{n}xT_{n-1}-\frac{n-1}{n}(x^2+y^2+z^2)T_{n-2},\ n\geq2,
    \end{aligned}
\end{equation}
with initial states $T_0=1,\ T_1=x$ and 
\begin{equation}
    \begin{aligned}
        R_n=\frac{2n+3}{n+2}xR_{n-1}-\frac{2n+2}{n+2}T_n-\frac{n+1}{n+2}(x^2+y^2+z^2)R_{n-2},\ n\geq2,
    \end{aligned}
\end{equation}
with initial states $R_0=-1,\ R_1=-3x.$
The coefficients $c_n$ depend solely on the system parameter $\mu$ and are given by
\begin{equation}
\begin{aligned}
    &c_n(\mu)=\frac{1}{{\gamma_i}^3}[(\pm1)^n\mu+(-1)^n\frac{(1-\mu){\gamma_i}^{n+1}}{(1\mp\gamma_i)^{n+1}}],\ \text{for}\ L_i,\ i=1,2; \\
    &c_n(\mu)=\frac{(-1)^n}{{\gamma_i}^3}[1-\mu+\frac{\mu{\gamma_i}^{n+1}}{(1+\gamma_i)^{n+1}}],\ \text{for}\ L_i,\ i=3.
\end{aligned}
\end{equation}
\section{Semi-analytical construction of bifurcated families in the ERTBP}\label{sec3}
In this section, we introduce a coupling coefficient and several bifurcation equations in the ERTBP. Based on different coupling directions, linear solutions are modified correspondingly. Initializing with these modified solutions, we develop a unified trigonometric series-based framework to iteratively construct a semi-analytical solution for describing the phase space near collinear libration points. 
By comprehensively considering all three cases of coupling constructions, the bifurcated orbits associated with the breaking of both $S_1$-type and $S_2$-type symmetries are approximated using modified high-order series expansions. 

The linearized equations associated with the ERTBP model (\ref{ERTBP3}) are given by
\begin{equation}\label{non-autonomous_linear_sol}
\begin{aligned}
    x^{\prime\prime}-2y^{\prime}-(1+2c_2)x=\sum\limits_{i\geq1}[(1+2c_2)x(-e)^i\cos^if],\\
    y^{\prime\prime}+2x^{\prime}+(c_2-1)y=\sum\limits_{i\geq1}[(1-c_2)y(-e)^i\cos^if],\\
    z^{\prime\prime}+c_2z=\sum\limits_{i\geq1}[(1-c_2)z(-e)^i\cos^if].
\end{aligned}
\end{equation}
Here, to obtain an explicit linear solution used in the initialization of the Lindstedt-Poincaré method, we start with the autonomous linear counterpart of (\ref{non-autonomous_linear_sol}), given by
\begin{equation}\label{autonomous_linear_sol}
    \begin{aligned}
        x^{\prime\prime}-2y^{\prime}-(1+2c_2)x=0,\\
        y^{\prime\prime}+2x^{\prime}+(c_2-1)y=0,\\
        z^{\prime\prime}+c_2z=0.
    \end{aligned}
\end{equation}
The first-order solution of (\ref{autonomous_linear_sol}) can be explicitly expressed as
\begin{equation}\label{non-modified_linear_sol}
    \begin{aligned}
        x(f)&=\alpha_1\cos{\theta_1}+\alpha_3\cos{\theta_3},\\
        y(f)&=\kappa_1\alpha_1\sin{\theta_1}+\sqrt{-1}\kappa_2\alpha_3\sin{\theta_3},\\
        z(f)&=\alpha_2\cos{\theta_2},
    \end{aligned}
\end{equation}
where $\alpha_1,\alpha_2\in \mathbb{R},\ \alpha_3\in {\sqrt{-1}\mathbb{R}}\cup{\mathbb{R}},\ \theta_1=\omega_0f+\varphi_1,\ \theta_2=\nu_0 f+\varphi_2, \ \text{and}\ \theta_3=\sqrt{-1}\lambda_0f+\varphi_3$, satisfying that
$$
    \begin{aligned}
        &\omega_0=\sqrt{\frac{2-c_2+\sqrt{9{c_2}^2-8c_2}}{2}},\ \nu_0=\sqrt{c_2},\ \lambda_0=\sqrt{\frac{c_2-2+\sqrt{9{c_2}^2-8c_2}}{2}};\\
        &\kappa_1=-\frac{{\omega_0}^2+2c_2+1}{2\omega_0},\ \kappa_2=-\frac{{\lambda_0}^2-2c_2-1}{2\lambda_0}.
    \end{aligned}
$$
Here, $\alpha_1$ and $\alpha_2$ represent the in-plane and out-of-plane amplitudes associated with the center part of the solution, respectively, while $\alpha_3$ corresponds to the amplitude associated with the hyperbolic part.
$\varphi_i$ ($i=1,2,3$) are the three corresponding initial phase angles. The coefficients $\kappa_1$ and $\kappa_2$ depend solely on the mass parameters of the RTBP. 

It is noted that the hyperbolic part in (\ref{non-modified_linear_sol}) is also expressed in a trigonometric form with complex amplitude and phase, where the value of $\alpha_3 \in {\sqrt{-1}\mathbb{R}}\cup{\mathbb{R}}$ is restricted to either the real axis or the imaginary axis. The solution with amplitude $\alpha_3$ lying on the real axis describes the motion of non-transit orbit, while the imaginary-valued $\alpha_3$ corresponds to the transit motion.

\textit{Remark 1.} Different from the works by Masdemont \cite{Masdemont} and Lei et al. \cite{Lei} where a hyperbolic exponential part is utilized to represent the motion in the hyperbolic direction, the trigonometric form applied here is helpful in analyzing bifurcations for solutions with hyperbolic components. Within this framework, the amplitude parameters $(\alpha_1,\alpha_2,\alpha_3)$ form bifurcation “indices” corresponding to each DOF of the ERTBP. By doing this, it is easier to implement a unified bifurcation analysis around not only periodic/quasi-periodic orbits but also transit/non-transit orbits. On the other hand, this trigonometric form contributes to a simpler formal expansion by taking advantage of the inherent symmetries of the ERTBP, as explained in the following subsections.

\subsection{Modification corresponding to the breaking of the $S_1$-type symmetry}
Along the family of planar Lyapunov orbits around collinear libration points, a pitchfork bifurcation occurs when the North-South symmetry (the $S_1$-type symmetry in (\ref{sy})) of the solution is broken. The resulting bifurcated classical orbits are known as halo orbits. To characterize this bifurcation, we consider a coupling effect from the motion in the $x$-direction to the $z$-direction in (\ref{ERTBP3}). In this case, the third equation in (\ref{ERTBP3}) is modified to:
\begin{equation}\label{modified_ERTBP}\tag{7.1}
    \begin{aligned}
z^{\prime\prime}+c_2z=&\sum\limits_{i\geq1}[-c_2 z(-e)^i\cos^if]
    +\sum\limits_{i\geq0}\{(-e)^i\cos^if[z\sum\limits_{n\geq2}c_{n+1}R_{n-1}(x,y,z)]\}+\eta\Delta x,
    \end{aligned}
\end{equation}
where $\eta$ is called the coupling coefficient and $\Delta$ represents the correction factor satisfying a bifurcation equation $\Delta=0$, which will be discussed later.
According to this introduced coupling effect, we obtain a modified liner equation, given by
\begin{equation}\label{x2z_modified_linear}
    \begin{aligned}
        x^{\prime\prime}-2y^{\prime}-(1+2c_2)x&=0,\\
        y^{\prime\prime}+2x^{\prime}+(c_2-1)y&=0,\\
        z^{\prime\prime}+c_2z&=\eta d_{0000}x.
    \end{aligned}
\end{equation}
By solving (\ref{x2z_modified_linear}), the modified linear solution is expressed as 
\begin{equation}\label{modified linear 1}
    \begin{aligned}
        x(f)&=\alpha_1\cos{\theta_1}+\alpha_3\cos{\theta_3},\\
        y(f)&=\kappa_1\alpha_1\sin{\theta_1}+\sqrt{-1}\kappa_2\alpha_3\sin{\theta_3},\\
        z(f)&=\alpha_2\cos{\theta_2}+\eta \alpha_1\cos{\theta_1}+\eta\kappa_3\alpha_3\cos{\theta_3},
    \end{aligned}
\end{equation}
where
$
\begin{aligned}
    \kappa_3=\frac{{\nu_0}^2-{\omega_0}^2}{{\nu_0}^2+{\lambda_0}^2},\ d_{0000}={\nu_0}^2-{\omega_0}^2.
\end{aligned}
$

When considering the perturbations of both nonlinear terms and orbital eccentricity, the high-order solution around collinear libration points in the ERTBP can be expressed as a formal expansion in powers of three amplitude parameters, and the orbital eccentricity:
\begin{equation}\label{formal sol}
    \begin{aligned}
    x(f)=&\sum x_{ijkm}^{stur}\cos{(s\theta_1+t\theta_2+u\theta_3+rf)}{\alpha_1}^i{\alpha_2}^j{\alpha_3}^k{e}^m,\\
    y(f)=&\sum y_{ijkm}^{stur}\sin{(s\theta_1+t\theta_2+u\theta_3+rf)}
    {\alpha_1}^i{\alpha_2}^j{\alpha_3}^k{e}^m,\\
    z(f)=&\sum z_{ijkm}^{stur}\cos{(s\theta_1+t\theta_2+u\theta_3+rf)}
    {\alpha_1}^i{\alpha_2}^j{\alpha_3}^k{e}^m,
    \end{aligned}
\end{equation}
where $\theta_1=\omega f+\varphi_1,\ \theta_2=\nu f+\varphi_2,\ \text{and}\ \theta_3=\sqrt{-1}\lambda f+\varphi_3$. Furthermore, the motion frequencies are not constant during the perturbation procedure and depend on the amplitudes $\alpha_i$ and the eccentricity $e$. Therefore, they are expressed as power series:
\begin{equation}\label{frequency}
    \begin{aligned}
        \omega=\sum\omega_{ijkm}{\alpha_1}^i{\alpha_2}^j{\alpha_3}^k{e}^m,\\
        \nu=\sum\nu_{ijkm}{\alpha_1}^i{\alpha_2}^j{\alpha_3}^k{e}^m,\\
        \lambda=\sum\lambda_{ijkm}{\alpha_1}^i{\alpha_2}^j{\alpha_3}^k{e}^m.
    \end{aligned}
\end{equation}
Likewise, the defined correction term $\Delta$ is expanded as $\Delta=\sum d_{ijkm}{\alpha_1}^i{\alpha_2}^j{\alpha_3}^k{e}^m$. 
To ensure that (\ref{formal sol}) is a valid solution of the original equation (\ref{ERTBP3}), the constraint condition $\eta\Delta=0$ must be satisfied.
Here, $\Delta=0$ establishes an implicit relationship between the amplitudes $\alpha_i$ ($i=1,2,3$), the orbital eccentricity $e$, and the coupling coefficient $\eta$. 
For any choice of the quartet $(\alpha_1,\alpha_2,\alpha_3,e)$, if there exist some $\eta\not=0$ satisfying the polynomial bifurcation equation $\Delta=0$, it indicates the occurrence of a bifurcation. 
Conversely, no bifurcation occurs if $\eta=0$. In this case, the high-order solution of planar/vertical Lyapunov orbits, Lissajous orbits and the corresponding transit/non-transit orbits can be derived. Since the iterative process is initialized with the modified linear solution (\ref{modified linear 1}), it satisfies that:
$$
    \begin{aligned}
        \omega_{0000}=\omega_0,\ \nu_{0000}=\nu_0,\ \lambda_{0000}=\lambda_0;\\
        x_{1000}^{1000}=x_{0010}^{0010}=1,\ y_{1000}^{1000}=\kappa_1,\ y_{0010}^{0010}=\sqrt{-1}\kappa_2;\\
        z_{1000}^{1000}=\eta,\ z_{0100}^{0100}=1,\ z_{0010}^{0010}
        =\eta\kappa_3.
    \end{aligned}
$$
\subsection{Modification corresponding to the breaking of the $S_2$-type symmetry}
Similar to halo families, it is known that associated with the breaking of the $S_2$-type symmetry in (\ref{sy}), the families of axial orbits bifurcate from both planar and vertical Lyapunov orbits. Each axial family consists of two branches related by the reflection across the $(x,y)$ plane. To illustrate this, a similar algorithm for the high-order solution of axial/quasi-axial orbits and their corresponding transit and non-transit orbits can be derived by introducing a coupling effect between the $z$-directional motion and the motion in the $y$-direction. Here, we can obtain the following two kinds of modified dynamical models.

For the first case, where the motion in the $y$-direction is coupled to the motion in the $z$-direction, the governing differential equations are obtained by correcting the third equation in (\ref{ERTBP3}) to
\begin{equation}\label{modified ERTBP2}\tag{7.2}
    \begin{aligned}
z^{\prime\prime}+c_2z=&\sum\limits_{i\geq1}[-c_2 z(-e)^i\cos^if]+\sum\limits_{i\geq0}\{(-e)^i\cos^if[z\sum\limits_{n\geq2}c_{n+1}R_{n-1}(x,y,z)]\}+\eta\Delta y,
    \end{aligned}
\end{equation}
with the modified linear solution used to initialize the perturbation procedure given by
\begin{equation}\label{y2z_modified_linear_sol}
    \begin{aligned}
       x(f)&=\alpha_1\cos{(\omega_0f+\varphi_1)}+\alpha_3\cos(\sqrt{-1}\lambda_0f+\varphi_3),\\
       y(f)&=\kappa_1\alpha_1\sin{(\omega_0f+\varphi_1)}+\sqrt{-1}\kappa_2\alpha_3\sin{(\sqrt{-1}\lambda_0f+\varphi_3)},\\
       z(f)&=\alpha_2\sin{(\nu_0f+\varphi_2)}+\eta\alpha_1\sin{(\omega_0f+\varphi_1)}+\sqrt{-1}\eta\kappa_3\alpha_3\sin{(\sqrt{-1}\lambda_0f+\varphi_3)},
    \end{aligned}
\end{equation}
where $d_{0000}=({\nu_0}^2-{\omega_0}^2)/\kappa_1,\ \kappa_3=\frac{\kappa_2}{\kappa_1}\frac{{\nu_0}^2-{\omega_0}^2}{{\nu_0}^2+{\lambda_0}^2}$.

For the second case, where the motion in the $z$-direction is coupled to the motion in the $y$-direction, the dynamical model is modified by reformulating the governing equation in the $y$-direction as 
\begin{equation}\label{modified ERTBP3}\tag{7.3}
    \begin{aligned}
    y^{\prime\prime}+2x^\prime+(c_2-1)y=\sum\limits_{i\geq1}[(1-c_2)y(-e)^i\cos^if]
    +\sum\limits_{i\geq0}\{(-e)^i\cos^if[y\sum\limits_{n\geq2}c_{n+1}R_{n-1}(x,y,z)]\}+\eta\Delta z.
    \end{aligned}
\end{equation}
In this case, the corresponding modified linear solution is expressed as
\begin{equation}
    \begin{aligned}
       x(f)&=\alpha_1\cos{(\omega_0f+\varphi_1)}+\eta\alpha_2\cos{(\nu_0f+\varphi_2)}+\alpha_3\cos(\sqrt{-1}\lambda_0f+\varphi_3),\\
       y(f)&=\kappa_1\alpha_1\sin{(\omega_0f+\varphi_1)}+\eta\kappa_3\alpha_2\sin{(\nu_0f+\varphi_2)}+\sqrt{-1}\kappa_2\alpha_3\sin{(\sqrt{-1}\lambda_0f+\varphi_3)},\\
       z(f)&=\alpha_2\sin{(\nu_0f+\varphi_2)},
    \end{aligned}
\end{equation}
with $d_{0000}=\frac{1}{2\nu_0}-\frac{\nu_0}{2},\ \kappa_3=-\frac{1}{2\nu_0}-\frac{3\nu_0}{2}.$ In both cases, the formal solution of the modified nonlinear dynamical model is expressed as 
\begin{equation}\label{formal solution type2}
    \begin{aligned}
    x(f)=&\sum x_{ijkm}^{stur}\cos{(s\theta_1+t\theta_2+u\theta_3+rf)}{\alpha_1}^i{\alpha_2}^j{\alpha_3}^k{e}^m,\\
    y(f)=&\sum y_{ijkm}^{stur}\sin{(s\theta_1+t\theta_2+u\theta_3+rf)}
    {\alpha_1}^i{\alpha_2}^j{\alpha_3}^k{e}^m,\\
    z(f)=&\sum z_{ijkm}^{stur}\sin{(s\theta_1+t\theta_2+u\theta_3+rf)}
    {\alpha_1}^i{\alpha_2}^j{\alpha_3}^k{e}^m,
    \end{aligned}
\end{equation}
while the expansions of the frequencies and the coupling correction term $\Delta$ retain the same formulations in (\ref{frequency}).

In contrast to the conventional perturbation techniques that depend on the 1:1 resonant modification of the in-plane and out-of-plane frequencies, our approach establishes a unified trigonometric series-based framework, where coupling-induced bifurcation mechanisms are used for systematically achieving symmetry breaking in non-bifurcated solutions.
The constructions of bifurcation equations according to different coupling effects are summarized with their associated symmetry breaking in Table \ref{table_sum}.

To obtain the semi-analytical solutions up to finite order $n$, coefficients associated with the formal expansions (\ref{formal sol}) and (\ref{formal solution type2}) are to be determined.
Relying on different modified linear solutions from preceding subsections as initial states, the computation of undetermined coefficients can be implemented iteratively using the Lindstedt-Poincaré method. The detailed computation associated with the breaking of the $S_1$-type symmetry is presented in Appendix, while other cases can be tackled in the same fashion.

\textit{Remark 2}. It is noticeable that unlike the method based on the 1:1 resonant bifurcation mechanisms, frequencies of both in-plane and out-of-plane motion are preserved in the coupling-bifurcation computation, enabling us to investigate the bifurcated halo/quasi-halo (axial/quasi-axial resp.) orbits in a unified framework. Moreover, taking advantage of the trigonometric expression of the hyperbolic part of the solution, the perturbation computation is also simplified compared to the conventional method. 

\textit{Remark 3}. It can be noticed that our analysis specifically addresses $S_1$- and $S_2$-type symmetry-breaking mechanisms, while $S_3$-type is absent.
In fact, $S_3$-type symmetry is a combination of the first two, meaning that its absence essentially corresponds to two sequential symmetry breaking.
Consequently, the orbits near collinear libration points may undergo two successive pitchfork bifurcations. 
From the global bifurcation diagram (Figure 3) of the CRTBP in \cite{Doedel}, this bifurcation corresponds to the $W_4$ and $W_5$ families of orbits.
The initial linear solutions for these orbits need to be obtained by incorporating corrections from two coupling-induced bifurcation equations. 
However, since the $W_4$ and $W_5$ families is significantly far from the collinear libration points, obtaining an accurate semi-analytical approximation of the $W_4$ and $W_5$ families using the local Lindstedt-Poincaré perturbation method remains challenging.

\renewcommand{\arraystretch}{1.8} 
\begin{table}
\caption{Classification of types of coupling directions and their corresponding symmetry breaking}\label{table_sum}
    \centering
    \begin{tabular}{ccc}
    \hline
     Type of symmetry breaking     & Coupling direction & Type of bifurcated orbit\\
     \hline
     $S_1$    & $x\to z$ & \makecell{Halo/quasi-halo orbits and their\\ corresponding transit/non-transit orbits}\\
     \hline
     $S_2$    & \makecell{$y\to z$ \\ $z\to y$} & \makecell{Axial/quasi-axial orbits and their \\corresponding transit/non-transit orbits}\\
     \hline
    \end{tabular}
\end{table}
\renewcommand{\arraystretch}{1}

\section{Results}\label{sec:results}
In this section, a detailed analysis for pitchfork bifurcations around the collinear libration points in the ERTBP is presented by tackling different parameterized bifurcation equations $\Delta(\eta,e,\alpha_1,\alpha_2,\alpha_3)=0$ quantitatively.
The emergence of non-zero solutions for $\eta$ induces bifurcated orbits, including periodic/quasi-periodic orbits, hyperbolic orbits, and transit/non-transit orbits, whose explicit critical conditions are also derived.

\subsection{Bifurcation associated with breaking of the \(S_1\)-type symmetry}
\subsubsection{Solvability of the third-order bifurcation equation}
\begin{figure}[!htb]
    \centering
    \begin{minipage}[t]{0.4\textwidth}
        \centering
        \includegraphics[width=\textwidth]{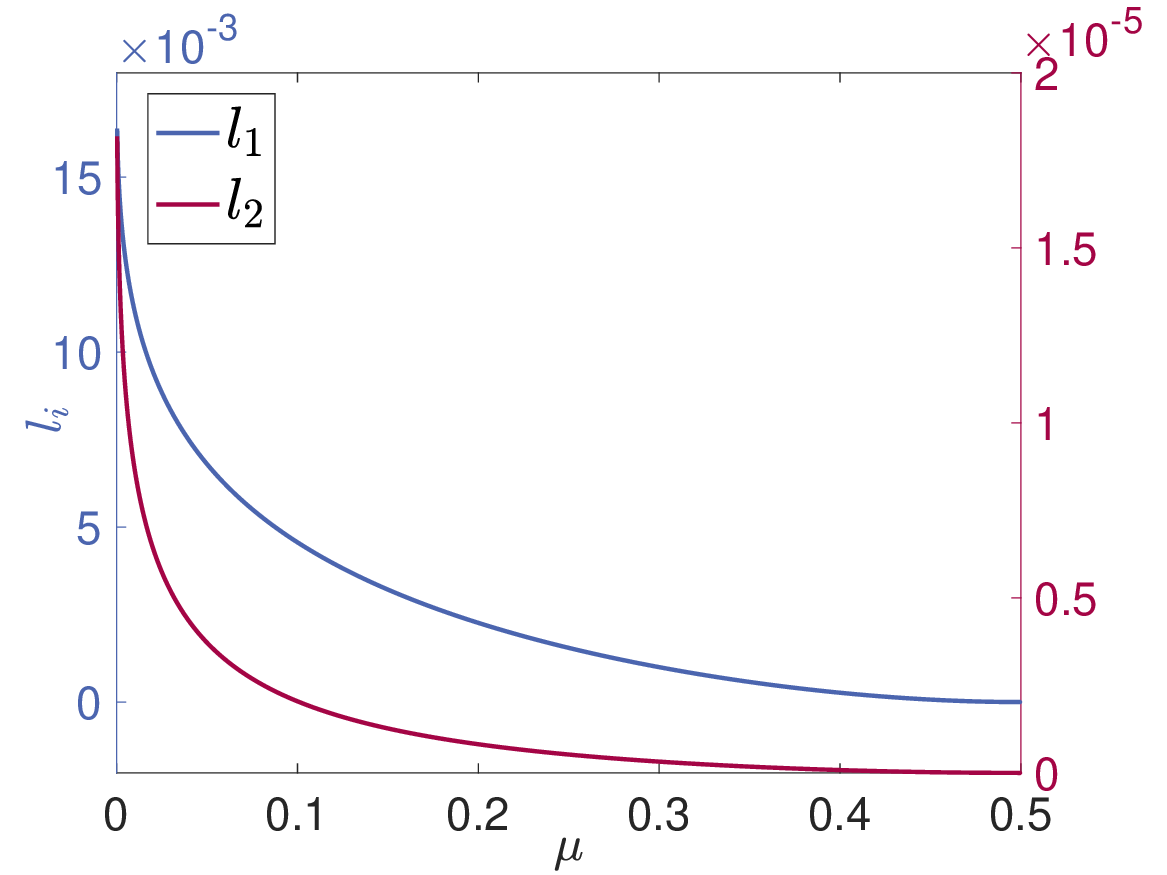}
    \end{minipage}
    \hspace{0.5cm}
    \begin{minipage}[t]{0.4\textwidth}
        \centering
        \includegraphics[width=\textwidth]{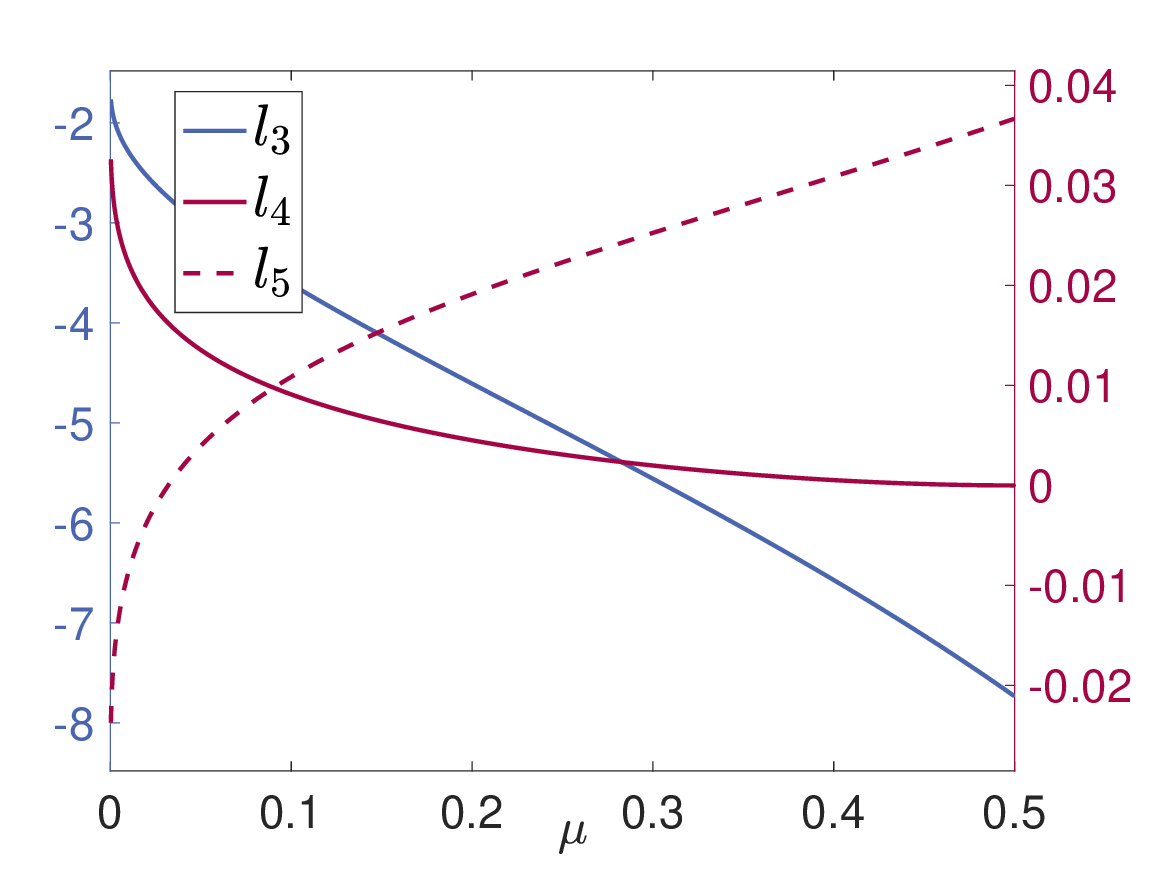}
    \end{minipage}
    
    \vspace{0.3cm} 
    
    \begin{minipage}[t]{0.4\textwidth}
        \centering
        \includegraphics[width=\textwidth]{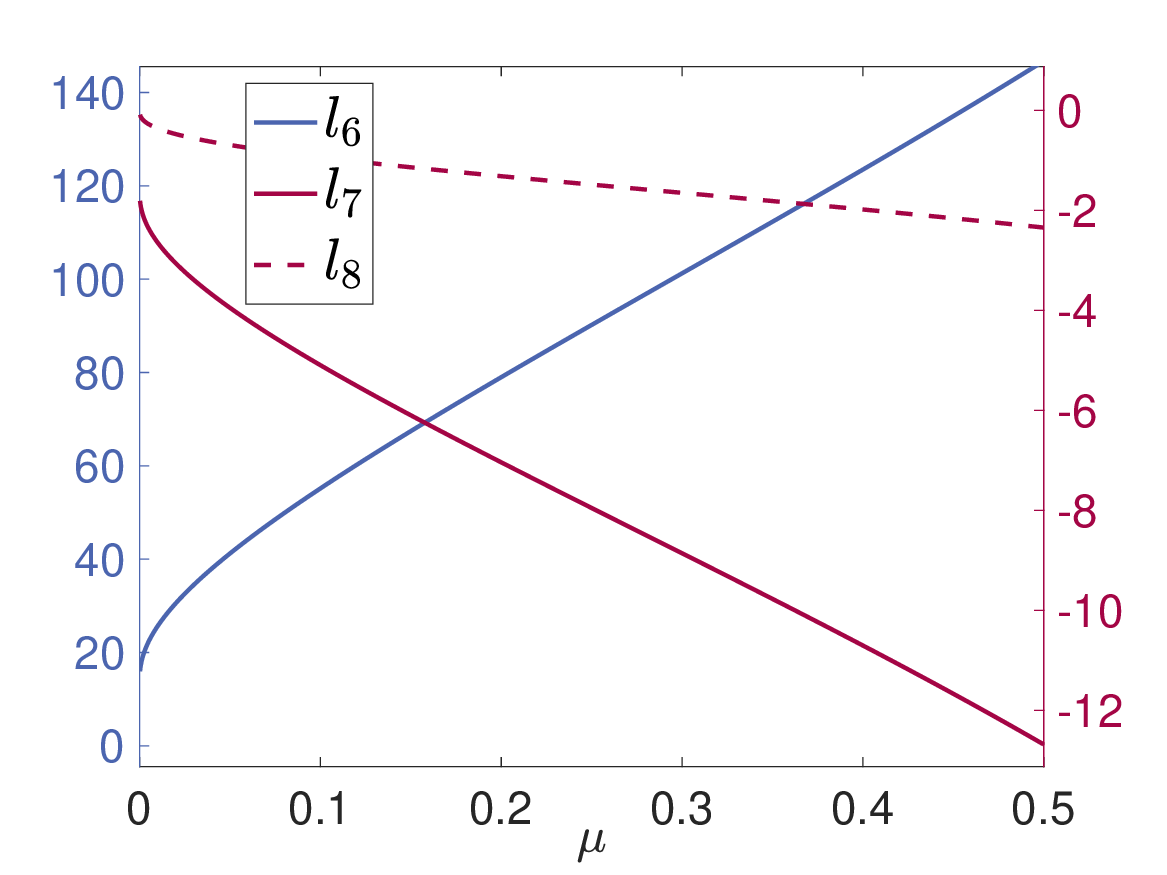}
    \end{minipage}
    \hspace{0.5cm}
    \begin{minipage}[t]{0.4\textwidth}
        \centering
        \includegraphics[width=\textwidth]{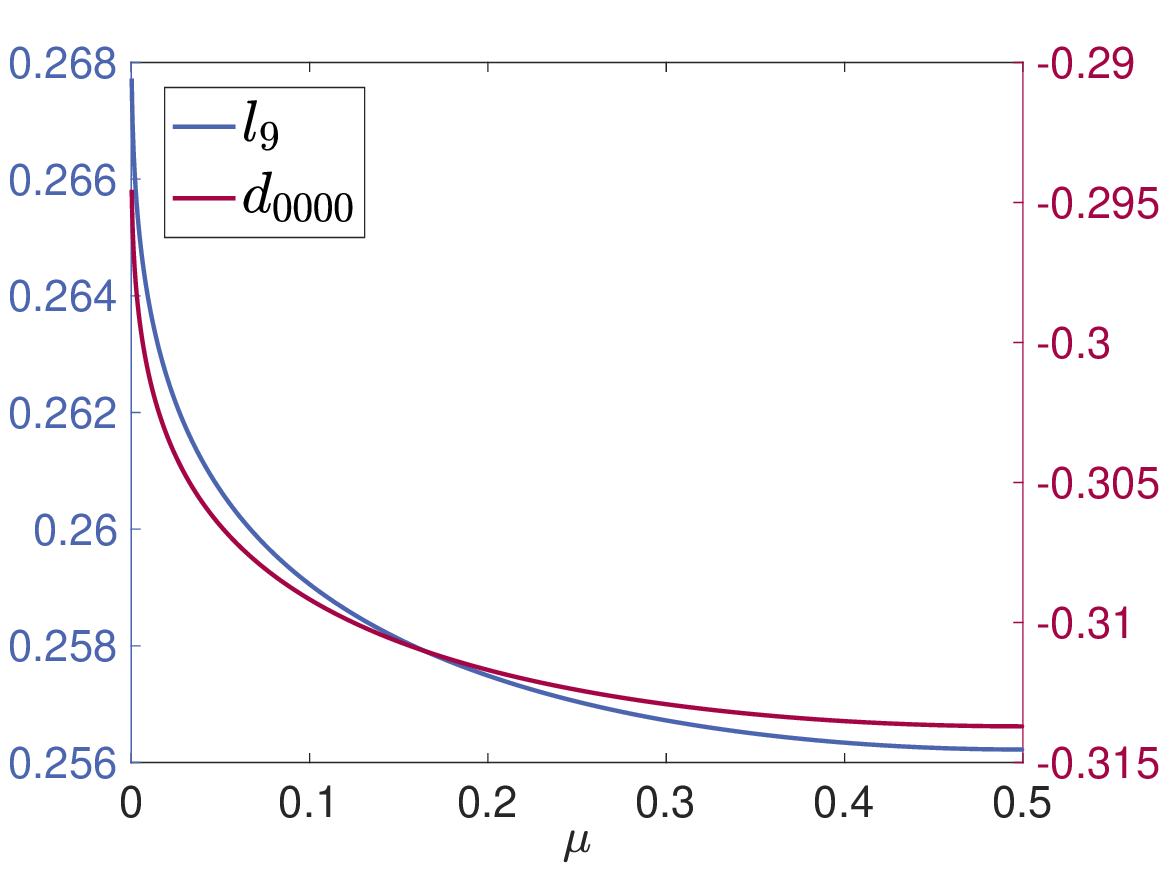}
    \end{minipage}
    \caption{Relationship between coefficients in the third-order bifurcation equation and the system parameter for $L_1$ in the ERTBP}\label{fig1}
\end{figure}

From the previous section, we know that the $S_1$-type symmetry breaking corresponding to a pitchfork bifurcation is induced by the coupling of motion in the x-direction to the z-direction. To analyze the specific bifurcation conditions, it is necessary to investigate solutions of the bifurcation equation $\Delta=0$. According to the linear equation for determining the coefficients $d_{ijkm}$ in (\ref{x2z_d_ijkm}), the bifurcation equation has non-zero solutions $\eta(\alpha_1,\alpha_2,\alpha_3,e)$ provided when the trigonometric series solution (\ref{formal sol}) is computed up to the third order, i.e.,
\begin{equation}\label{bifurcation equation1}
\begin{aligned}
        \Delta&=d_{0000}+d_{2000}{\alpha_1}^2+d_{0200}{\alpha_2}^2+d_{0020}{\alpha_3}^2+d_{0002}e^2\\
        &=a\eta^4+b\eta^2+c=0,
\end{aligned}
\end{equation}
where $a=l_1{\alpha_1}^2+l_2{\alpha_3}^2,b=l_3{\alpha_1}^2+l_4{\alpha_2}^2+l_5{\alpha_3}^2,c=l_6{\alpha_1}^2+l_7{\alpha_2}^2+l_8{\alpha_3}^2+l_9{e}^2+({\nu_0}^2-{\omega_0}^2)$. The coefficients $l_i$ ($i=1,2,...,9$) depend solely on the system parameter $\mu$ and their quantitative relationships are illustrated in Fig.~\ref{fig1}. 
Here, the relatively negligible impact of the small orbital eccentricity on the bifurcation equation illustrates some significant similarities shared by the non-autonomous ERTBP and its approximated circular model from quantitative perspectives.
By treating the bifurcation equation (\ref{bifurcation equation1}) as a quadratic equation in $\eta^2$, a detailed bifurcation analysis can be implemented as follows.

As previously discussed, the complex-valued amplitude $\alpha_3$ which corresponds to the hyperbolic motion comprises two branches while one refers to non-transit trajectories (i.e. $\alpha_3\in \mathbb{R}$) and the other branch describes the transit trajectories (i.e. $\alpha_3\in \sqrt{-1}\mathbb{R}$). For the sake of simplicity, these two cases will be analyzed separately.\\
\textbf{Case 1.1:} $\alpha_3\in \sqrt{-1}\mathbb{R},c=0,-\frac{b}{a}>0$.\\
Replace the imaginary-valued amplitude $\alpha_3$ with $\tilde{\alpha}_3:=\alpha_3/\sqrt{-1}$. Hereafter, $\alpha_3$ will denote its 
real-valued counterpart. The critical surface $c=0$ is then formulated as 
\begin{equation}\label{case1.1}
    c=l_6{\alpha_1}^2+l_7{\alpha_2}^2-l_8{\alpha_3}^2+l_9{e}^2+({\nu_0}^2-{\omega_0}^2)=0.
\end{equation}
For any $\mu\in (0,0.5)$, the coefficients satisfy $l_6>0;\ l_7,l_8<0$. Now, (\ref{case1.1}) defines a one-sheet hyperboloid in the ($\alpha_1,\alpha_2,\alpha_3$) coordinate system for any orbital eccentricity $0<e<1$. On the hyperboloid surface with condition $-b/a>0$, there exist two distinct non-zero solutions given by $\eta=\pm\sqrt{-\frac{b}{a}}$. On the right side of the hyperboloid, four feasible $\eta$ are obtained as 
\begin{equation}
    \eta=\pm\sqrt{\frac{-b\pm\sqrt{b^2-4ac}}{2a}},
\end{equation}
which solve the quadratic equation, while on the left side of the critical surface, the negativity of the coefficient $c$ restricts the solutions to only two feasible values: $\eta=\pm\sqrt{\frac{-b+\sqrt{b^2-4ac}}{2a}}$. \\
\textbf{Case 1.2:} $\alpha_3\in \sqrt{-1}\mathbb{R},a=0,-\frac{c}{b}>0.$\\
In this case, the critical surface is represented as 
\begin{equation}
    a=l_1{\alpha_1}^2-l_2{\alpha_3}^2=0.
\end{equation}
On the critical surface, two distinct solutions $\eta=\pm\sqrt{-\frac{c}{b}}$ emerge. On the left part, there exist four feasible solutions, whereas on the other part, two solutions corresponding to (\ref{bifurcation equation1}) are given by $\eta=\pm\sqrt{\frac{-b+\sqrt{b^2-4ac}}{2a}}.$

\begin{figure}[!htb]
    \centering
    \subfigure[]{
    \includegraphics[width=0.4\textwidth]{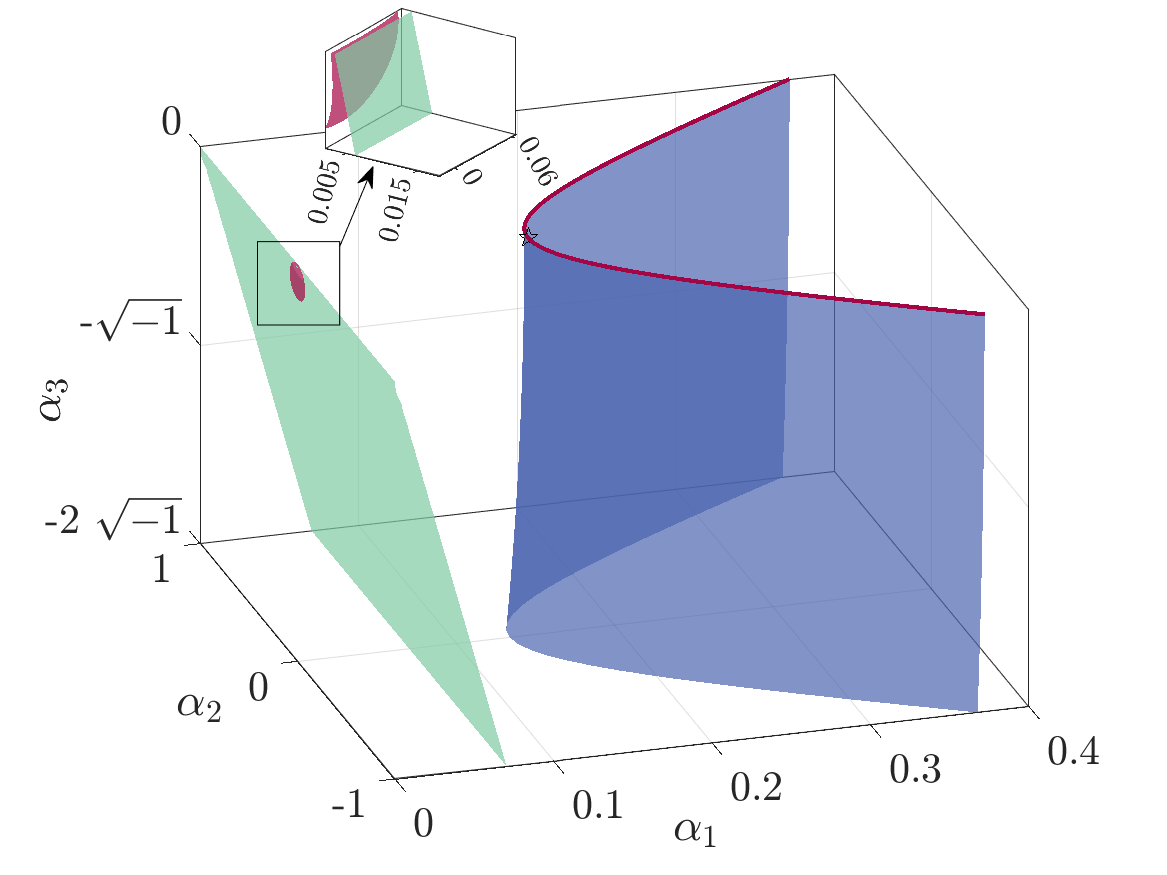}
    }
    \subfigure[]{
    \includegraphics[width=0.4\textwidth]{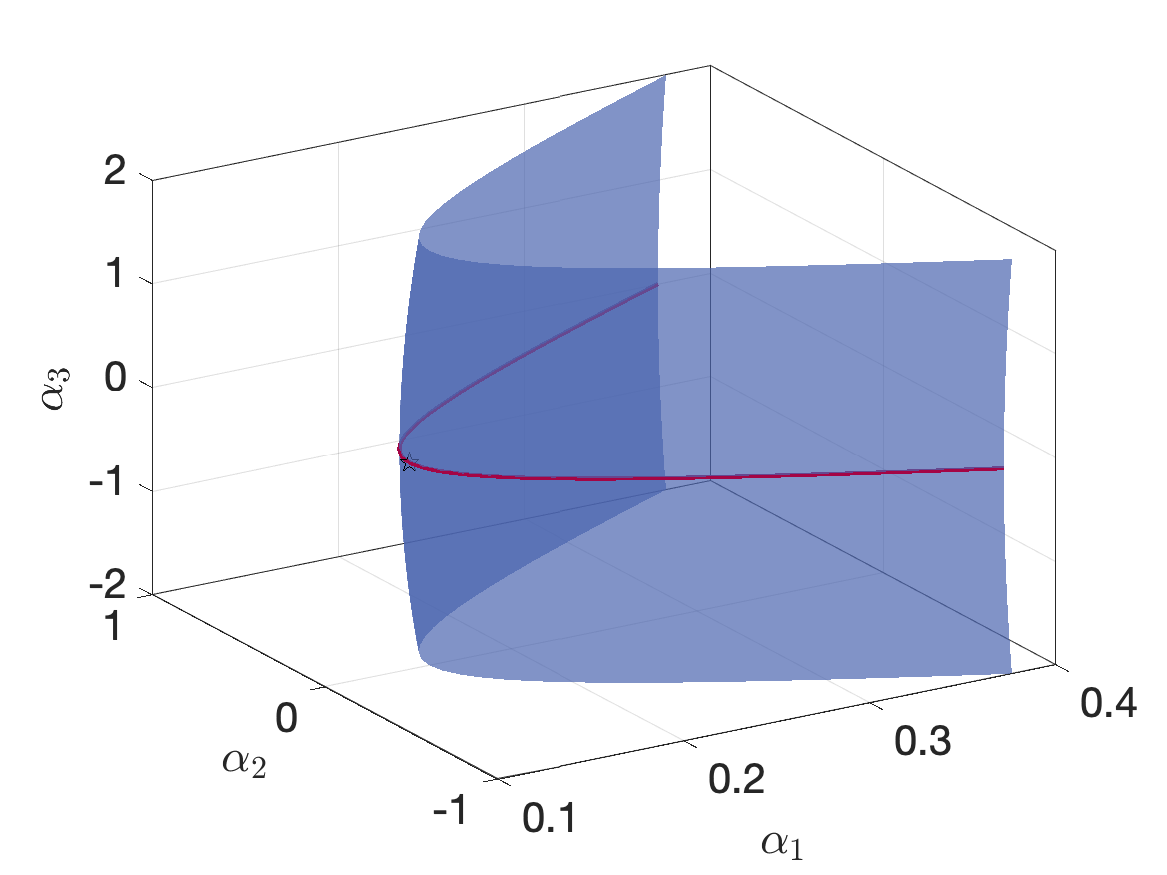}
    }
    \caption{Distribution of feasible solutions of the third-order bifurcation equation in the Sun-Earth system. (a): Critical surfaces corresponding to transit orbits. (b): Critical surfaces corresponding to non-transit orbits.}
    \label{fig2}
\end{figure}

It is noteworthy that both \textbf{Case 1.1} and \textbf{Case 1.2} suggest the occurrence of some potential bifurcations in the hyperbolic part of the solution and, more generally, in transit and non-transit orbits. Unrevealed from the resonant bifurcation mechanisms, these newly identified bifurcated orbits in the ERTBP can now be systematically characterized through the parameterized bifurcation equation. Detailed descriptions of these orbits are provided in the subsequent subsections.\\
\textbf{Case 1.3:} $\alpha_3\in\sqrt{-1}\mathbb{R},b^2-4ac=0,-\frac{b}{a}>0.$\\
On the complicated critical surface defined by
\begin{equation}
    (l_3{\alpha_1}^2+l_4{\alpha_2}^2-l_5{\alpha_3}^2)^2-4(l_1{\alpha_1}^2-l_2{\alpha_3}^2)[l_6{\alpha_1}^2+l_7{\alpha_2}^2-l_8{\alpha_3}^2+l_9{e}^2+({\nu_0}^2-{\omega_0}^2)]\\
    =0,
\end{equation}
$\eta=\pm\sqrt{-\frac{b}{2a}}$ are the two feasible solutions corresponding to the bifurcation equation. Inside the surface, the condition $b^2-4ac<0$ holds, indicating no bifurcation occurs in this region. Fig.~\ref{fig2}(a) illustrates these three critical cases as introduced above in the Sun-Earth system where the system parameter $\mu=3.040423398444176$e-6 and the orbital eccentricity $e=0.01671022$. Four feasible solutions exist in the region to the right of the blue surface. Two solutions exist between the blue and green surfaces, and four solutions satisfy the bifurcation equation to the left of the green surface. Inside the red surface, no feasible solution exists. 

Similarly, the solvability of (\ref{bifurcation equation1}) can be analyzed for amplitudes associated with non-transit orbits (i.e. $\alpha_3\in \mathbb{R}$) in a manner analogous to the transit case. Specifically, the first critical surface is comparably defined as $\{(\alpha_1,\alpha_2,\alpha_3,e)\in \mathbb{R}^3\times(0,1))|c=0,-\frac{b}{a}>0\}$. In this case, the surface $c=0$ describes a two-sheet hyperboloid in the ($\alpha_1,\alpha_2,\alpha_3$) coordinate system, as shown in Fig.~\ref{fig2}(b). Similar to \textbf{Case 1.1}, two additional feasible solutions, $\eta=\pm\sqrt{\frac{-b-\sqrt{b^2-4ac}}{2a}}$ emerge on the right side of the critical surface. For the non-transit case, the coefficient $a=l_1{\alpha_1}^2+l_2{\alpha_3}^2$ is positive everywhere except at the origin, ensuring the degenerate case $\{(\alpha_1,\alpha_2,\alpha_3,e)\in \mathbb{R}^3\times(0,1))|a=0,-c/b>0\}$ does not exist.

\subsubsection{Bifurcation analysis restricted to center manifolds}
\begin{figure}[!htb]
    \centering
    \subfigure[]{
    \includegraphics[width=0.4\textwidth]{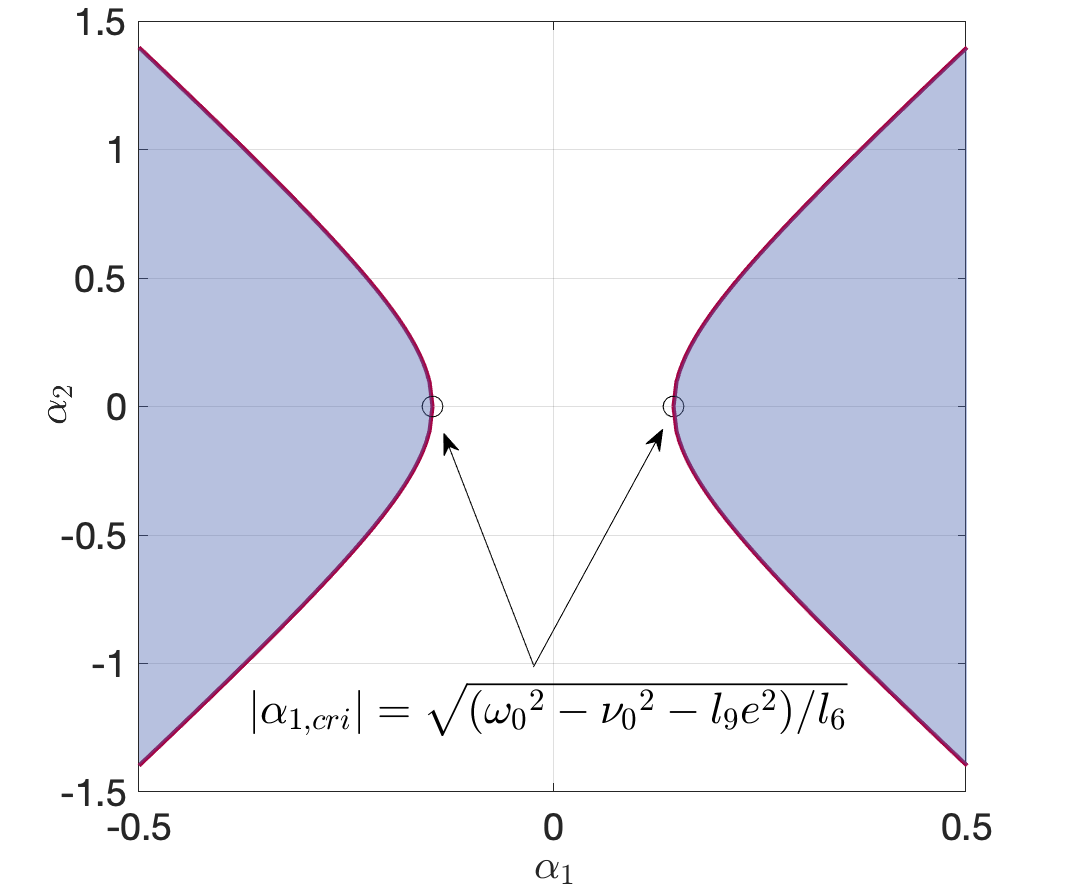}
    }
    \subfigure[]{
    \includegraphics[width=0.4\textwidth]{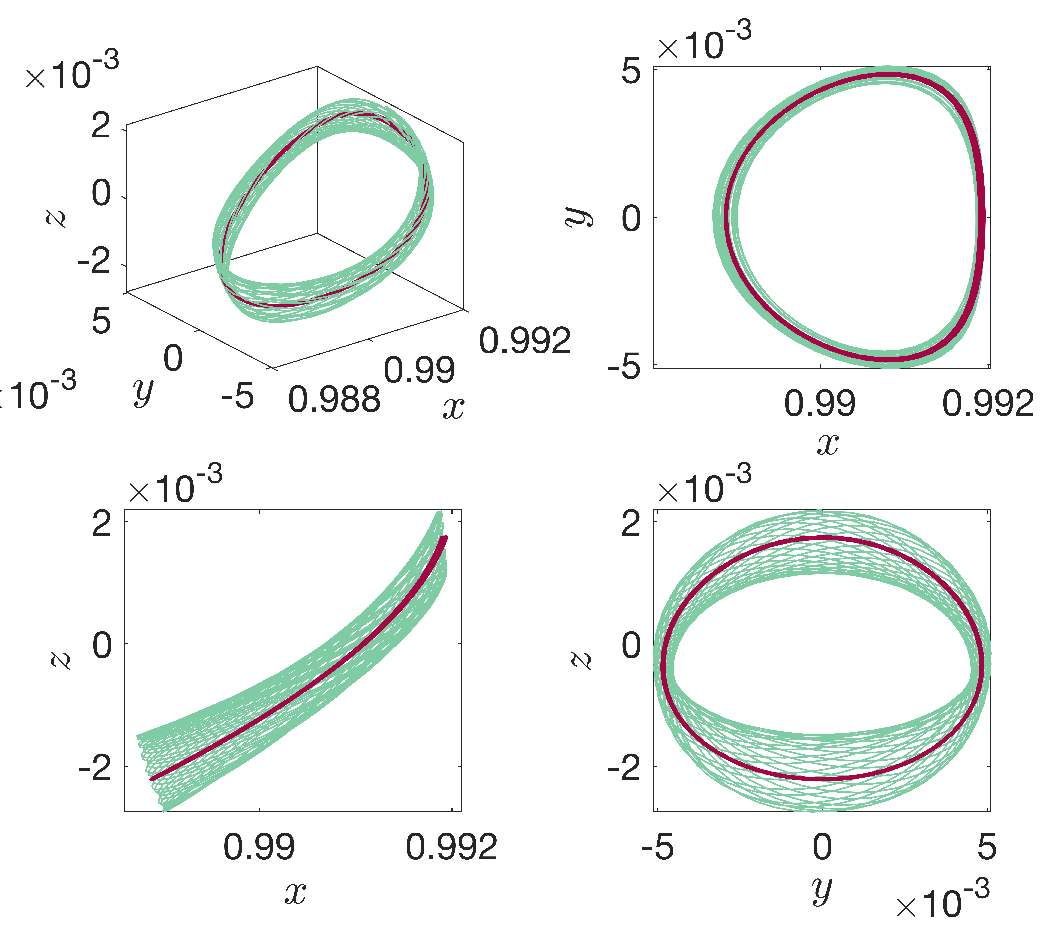}
    }
    \caption{(a): Bifurcation diagram of the third-order series solution of orbits in center manifolds associated with breaking of the $S_1$-type symmetry. (b): Red: Halo orbit with $\alpha_1=0.15$. Green: Quasi-halo orbit with $\alpha_1=0.15,\ \alpha_2 = 0.04$.}
    \label{fig3}
\end{figure}

By setting $\alpha_3=0$ in (\ref{bifurcation equation1}), the reduced bifurcation equation can be expressed as 
\begin{equation}\label{bifurcation equation center}
    \Delta=\hat{a}\eta^4+\hat{b}\eta^2+\hat{c}=0,
\end{equation}
where $\hat{a}=l_1{\alpha_1}^2,\hat{b}=l_3{\alpha_1}^2+l_4{\alpha_2}^2,\ \text{and}\ \hat{c}=l_6{\alpha_1}^2+l_7{\alpha_2}^2+l_9{e}^2+({\nu_0}^2-{\omega_0}^2)$. In this case, the bifurcation curve in the $(\alpha_1,\alpha_2)$ plane is defined by ${\hat{c}=0}$. For any system parameter $\mu\in (0,0.5)$ and a small orbital eccentricity $e$, the bifurcation curve
\begin{equation}
    l_6{\alpha_1}^2+l_7{\alpha_2}^2={\omega_0}^2-{\nu_0}^2-l_9{e}^2
\end{equation}
describes a hyperbola, as illustrated in Fig.~\ref{fig3}(a), where two critical points ($\pm\sqrt{({\omega_0}^2-{\nu_0}^2-l_9{e}^2)/l_6},0$) lie on the $\alpha_1$-axis.
When $\eta=0$ and $\Delta\not=0$, the series expansions (\ref{formal sol}) provide approximations to non-bifurcated orbits, encompassing Lyapunov orbits and Lissajous orbits.    
However, as $\alpha_1$ and $\alpha_2$ increase, non-zero solutions $\eta$ to the bifurcation equation $\Delta=0$ may emerge.
These bifurcated solutions describe halo and quasi-halo orbits. Specifically, when $\alpha_2=0$, the series solutions with $\eta=0$ represent the family of planar Lyapunov orbits. For $|\alpha_1|>\sqrt{({\omega_0}^2-{\nu_0}^2-l_9{e}^2)/l_6}$, (\ref{bifurcation equation center}) yields a pair of non-zero solutions given by
\begin{equation}
\eta=\pm\sqrt{\frac{-l_3{\alpha_1}^2-\sqrt{{l_3}^2{\alpha_1}^4-4l_1{\alpha_1}^2(l_6{\alpha_1}^2+l_9{e}^2+({\nu_0}^2-{\omega_0}^2))}}{2l_1{\alpha_1}^2}},
\end{equation}
corresponding to the northern halo orbits and southern halo orbits respectively. When $\alpha_2\not=0$, Lissajous orbits are defined with $\eta=0$, while quasi-halo orbits bifurcate from Lissajous orbits if $\eta\not=0$ satisfies the bifurcation equation $\Delta=0$.
By using the Lindstedt-Poincaré method introduced in the preceding sections, here, the semi-analytical solution is computed up to the 7th order. Halo/quasi-halo orbits bifurcated from planar Lyapunov/Lissajous orbits around $L_1$ in the Sun-Earth ERTBP are shown in Fig.~\ref{fig3}(b).
\subsubsection{Bifurcation analysis of hyperbolic orbits}
\begin{figure}[!htb]
    \centering
    \subfigure[]{
    \includegraphics[width=0.4\textwidth]{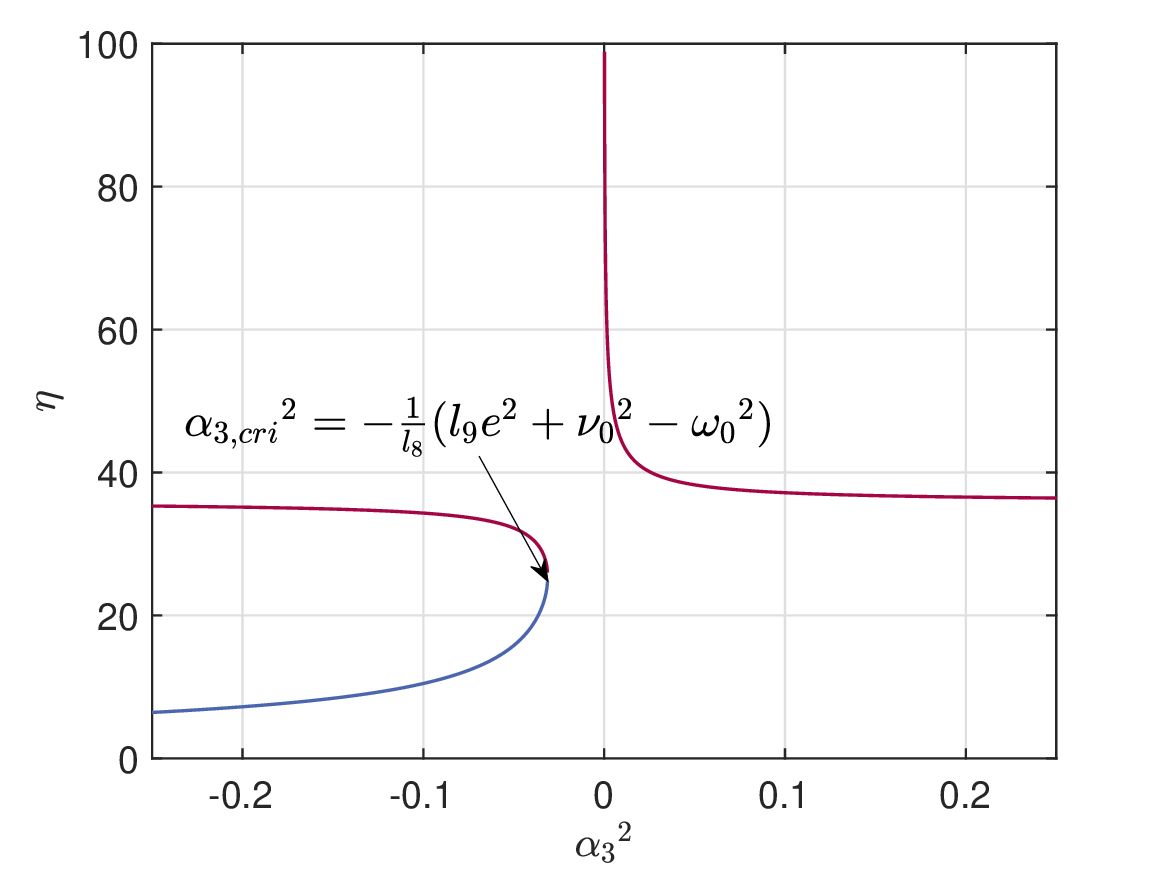}
    }
    \subfigure[]{
    \includegraphics[width=0.4\textwidth]{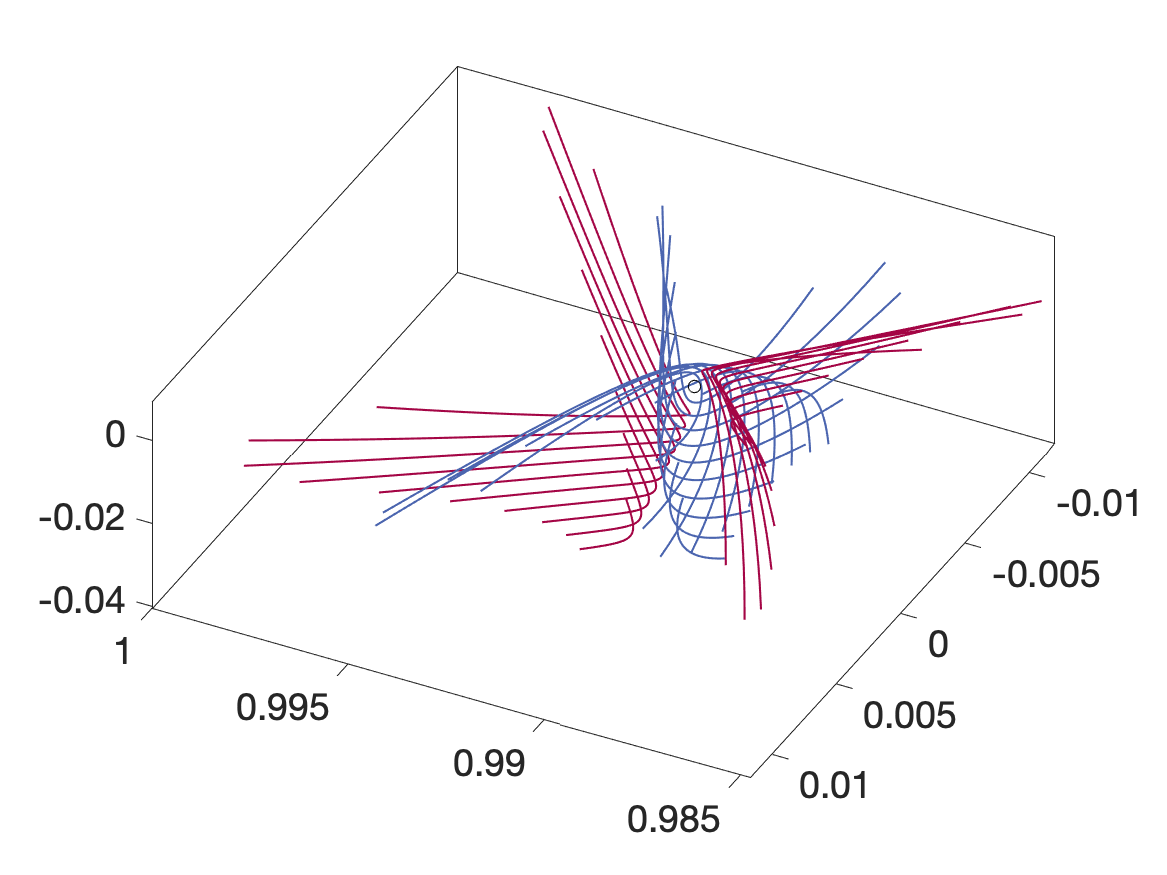}
    }
    \caption{(a): Positive solutions of the bifurcation equation with $\alpha_1=\alpha_2=0$ corresponding to $L_1$ in the Sun-Earth system. (b): Hyperbolic orbits bifurcated from planar hyperbolic orbits.}
    \label{fig4}
\end{figure}

By letting the amplitudes corresponding to center part of (\ref{formal sol}), $\alpha_1,\alpha_2$ be zero, the bifurcation equation is simplified to
\begin{equation}\label{bifurcation equation hyperbolic}
    \Delta=l_2{\alpha_3}^2\eta^4+l_5{\alpha_3}^2\eta^2+l_8{\alpha_3}^2+l_9{e}^2+({\nu_0}^2-{\omega_0}^2)=0.
\end{equation}
When $\eta=0$, $\Delta\not=0$, no bifurcation occurs. In this case, (\ref{formal sol}) describes planar hyperbolic orbits. Bifurcated hyperbolic orbits emerge from these planar hyperbolic orbits when $\eta$ is chosen as a non-zero solution to (\ref{bifurcation equation hyperbolic}). Specifically, for ${\alpha_3}^2>0$, there exist a pair of real-valued $\eta$:
\begin{equation}
\eta=\pm\sqrt{\frac{-l_5{\alpha_3}^2+\sqrt{{l_5}^2{\alpha_3}^4-4 l_2{\alpha_3}^2[l_8{\alpha_3}^2+l_9{e}^2+({\nu_0}^2-{\omega_0}^2)]}}{2 l_2{\alpha_3}^2}},
\end{equation}
which solves the reduced bifurcation equation. Conversely, for ${\alpha_3}^2<0$ and $l_8{\alpha_3}^2+l_9{e}^2+({\nu_0}^2-{\omega_0}^2)>0$, we can obtain up to four feasible solutions for (\ref{bifurcation equation hyperbolic}), as seen in Fig.~\ref{fig4}(a). Fig.~\ref{fig4}(b) presents the bifurcated hyperbolic orbits that evolve beyond the $(x,y)$ plane.
\subsubsection{Bifurcation analysis of transit/non-transit orbits}
\begin{figure}[!htb]
    \centering
    \subfigure[]{
    \includegraphics[width=0.4\textwidth]{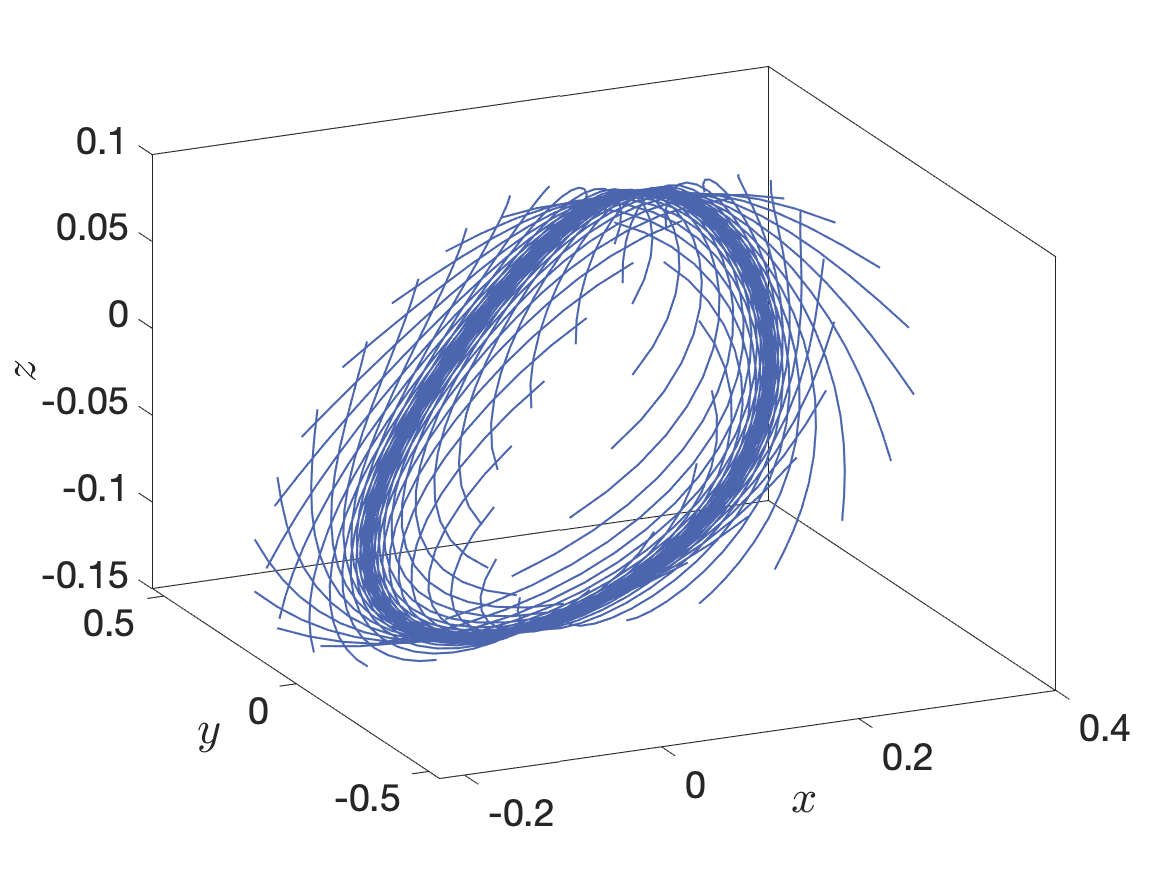}
    }
    \subfigure[]{
    \includegraphics[width=0.4\textwidth]{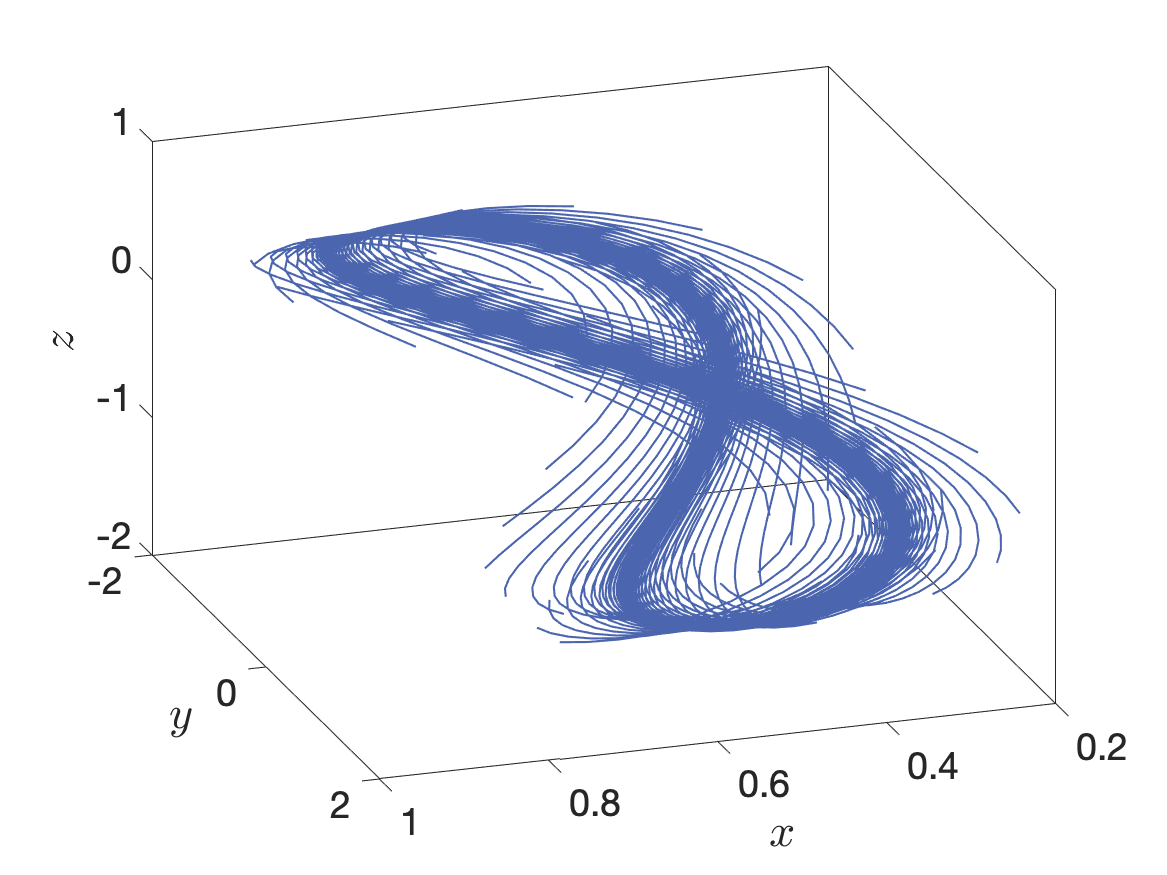}
    }
    \caption{Two branches of bifurcated transit orbits in the Sun-Earth system. (a) Transit orbits with $\alpha_1=0.15,\ \alpha_3=\pm0.005\sqrt{-1}$,\ $\eta=0.7085$. (b) Transit orbits with $\alpha_1=0.15,\ \alpha_3=\pm0.005\sqrt{-1}$,\ $\eta=9.4484$.}
    \label{fig5}
\end{figure}

As previously discussed, when the amplitude $\alpha_3$ associated with hyperbolic manifolds satisfies $\alpha_3\in \sqrt{-1}\mathbb{R}$, the series solution (\ref{formal sol}) characterizes transit orbits in the ERTBP. Spacecraft can transit from one side of the collinear libration points to the other along these orbits. It is obtained from the bifurcation equation (\ref{bifurcation equation hyperbolic}) that transit orbits also undergo bifurcations.
Here, depending on different choices of feasible solutions $\eta$ to $\Delta=0$, we present two families of bifurcated transit orbits with different dynamical behaviors. 
Now, consider the case where $\alpha_1\not=0,\alpha_2=0$. When $\eta\not=0$, It is observed that transit orbits bifurcate from planar orbits. Specifically, for relatively large values of $\eta$, the dynamics of the bifurcated transit orbits are primarily governed by motion in the hyperbolic direction. These orbits exhibit relatively rapid escape from the $(x,y)$ plane. When $\eta$ is selected from pairs with small absolute values $|\eta|$, an additional branch of bifurcated transit orbits (known as transit orbits of halo orbits in the ERTBP) emerges. In this scenario, the motion in the center part dominates, resulting in a slow variation of the $z$-axis. Fig.~\ref{fig5} illustrates the two distinct families of transit orbits bifurcated from planar orbits, highlighting the contrasting dynamical behaviors associated with different values of $\eta$.
\subsection{Bifurcation associated with the breaking of the $S_2$-type symmetry}
\subsubsection{Case of coupling the motion in the $z$-direction to the motion in the $y$-direction.}
\begin{figure}[!htb]
    \centering
    \subfigure[]{
    \includegraphics[width=0.4\textwidth]{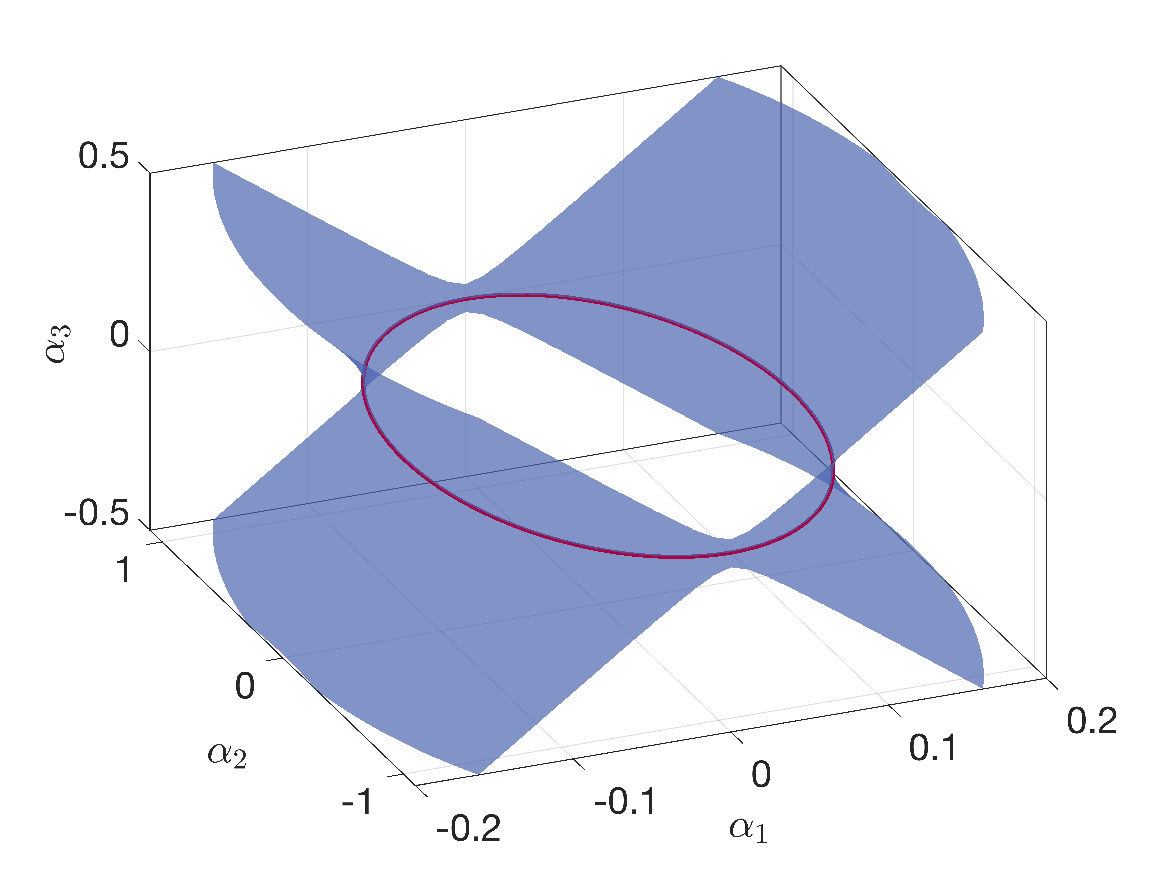}
    }
    \subfigure[]{
    \includegraphics[width=0.4\textwidth]{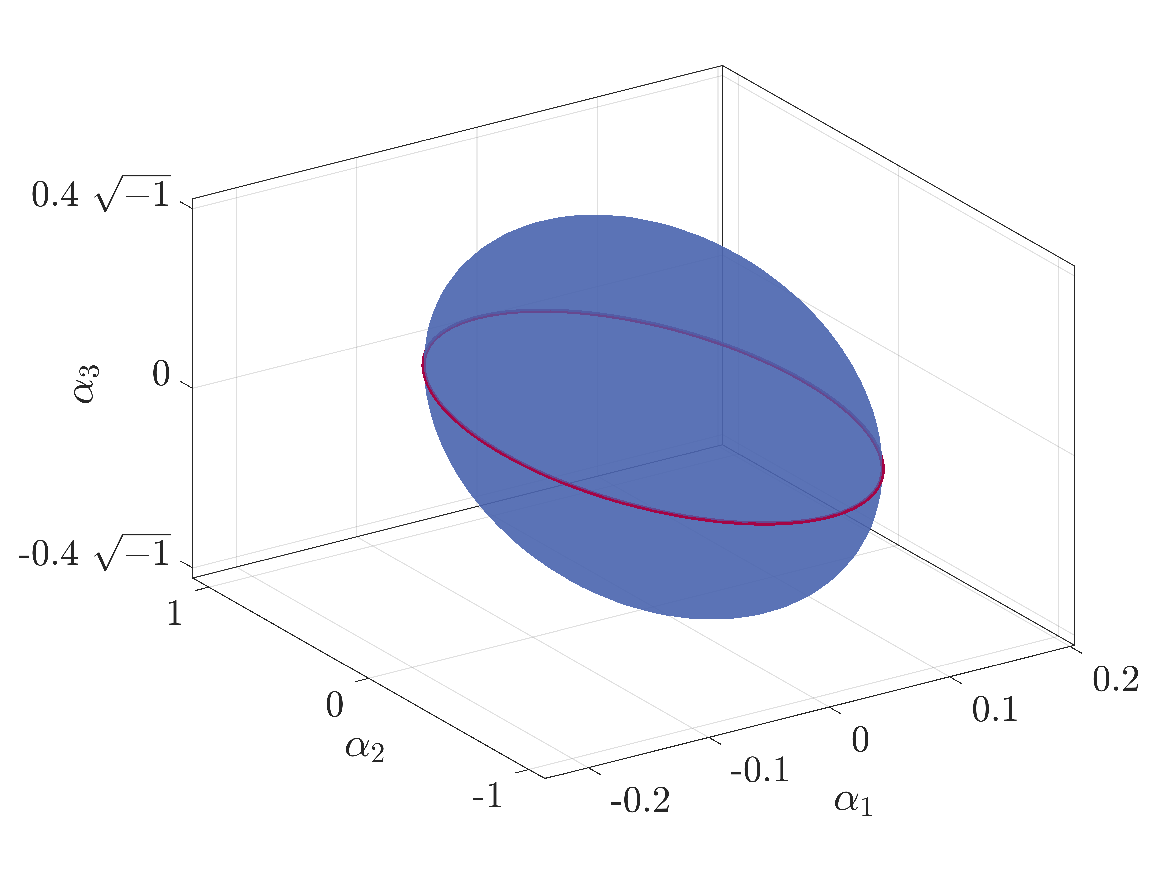}
    }
    \caption{Critical surfaces determined by the bifurcation equation for the case of considering a coupling effect from the motion in the $z$-direction to the $y$-direction. (a) Critical surface with $\alpha_3$ associated with non-transit orbits. (b) Critical surface with $\alpha_3$ associated with transit orbits.}
    \label{fig6}
\end{figure}

\begin{figure}[!htb]
    \centering
    \subfigure[]{
    \includegraphics[width=0.4\textwidth]{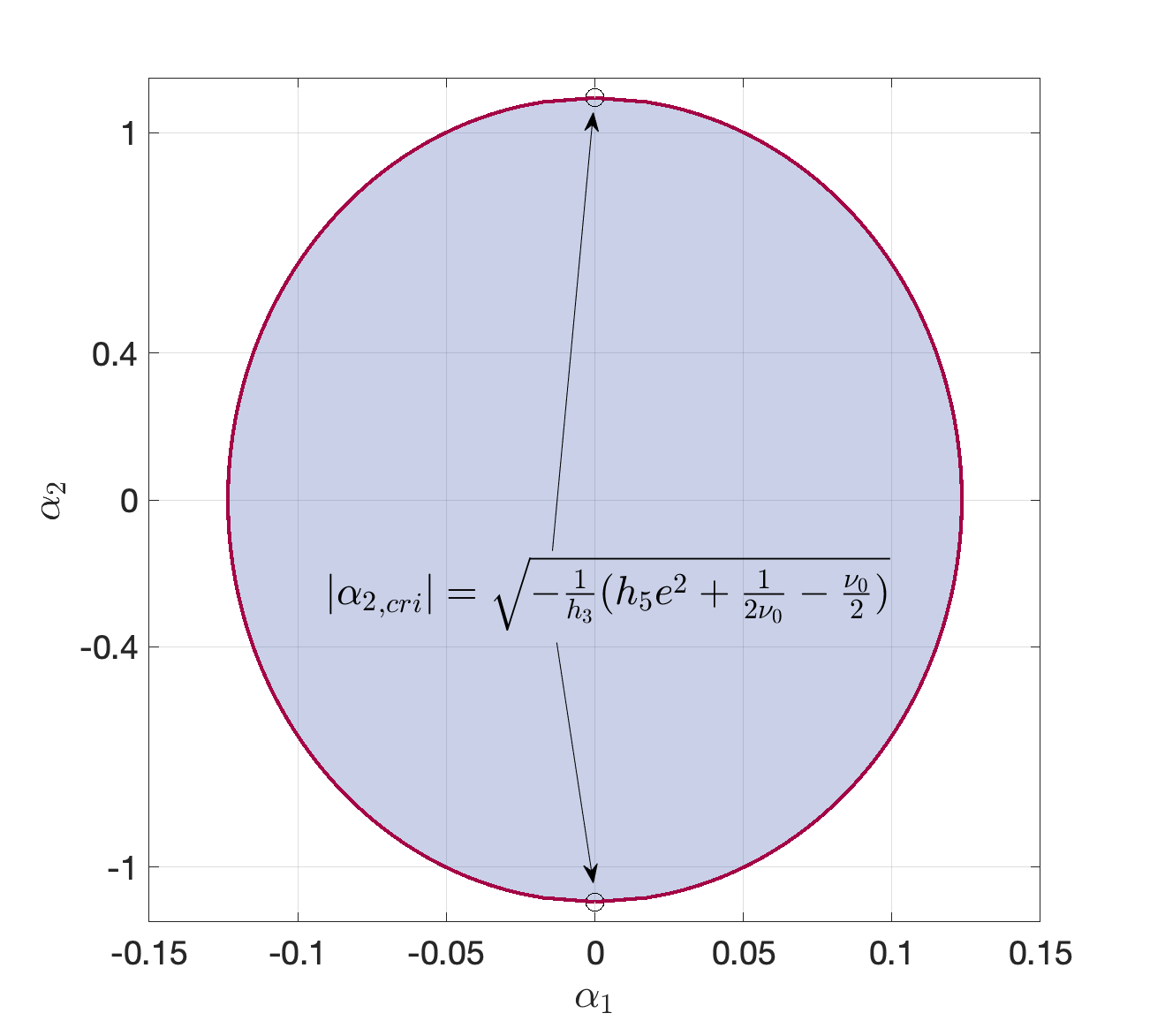}
    }
    \subfigure[]{
    \includegraphics[width=0.4\textwidth]{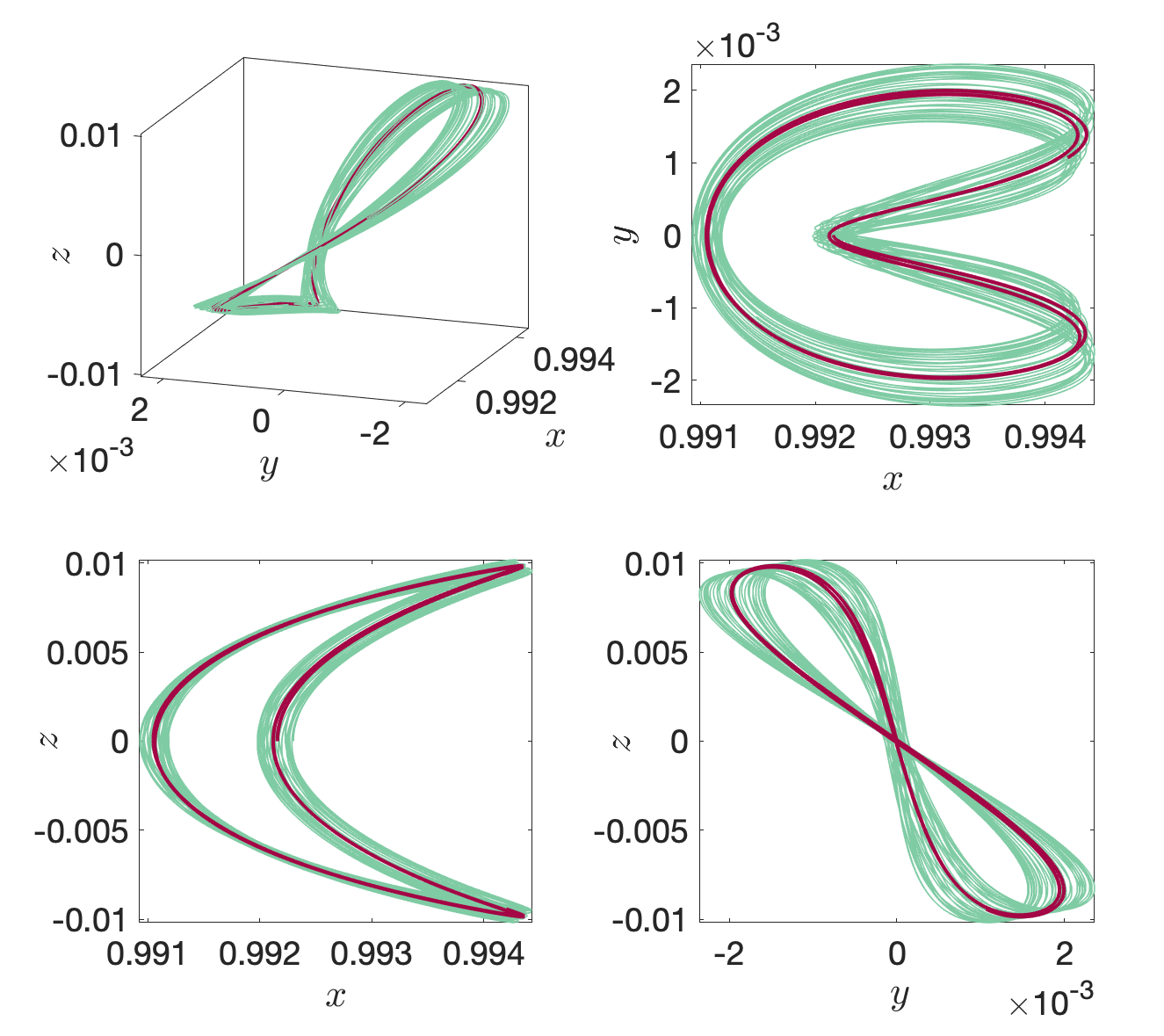}
    }
    \caption{(a): Bifurcation diagram of the third-order series solution of orbits in center manifolds in the case of coupling the motion in the $z$-direction to the motion in the $y$-direction. (b): Red: Axial orbit with $\alpha_1=0,\ \alpha_2=1$. Green: Quasi-halo orbit with $\alpha_1=0.005,\ \alpha_2 = 1$.}
    \label{fig7}
\end{figure}

In this case, the bifurcated solution first appears in the third-order series solution. The bifurcation equation is given by
\begin{equation}\label{bifurcation equation z2y}
    \Delta=h_1{\alpha_2}^2\eta^2+h_2{\alpha_1}^2+h_3{\alpha_2}^2+h_4{\alpha_3}^2+h_5e^2+\frac{1}{2\nu_0}-\frac{\nu_0}{2}=0,
\end{equation}
where $h_i\ (1\leq i\leq 5)$ depends solely on the system parameter $\mu$. 
In the Sun-Earth system, the coefficients satisfy $h_2,h_3>0,\ h_4<0$. The critical surface determined by the bifurcation equation (\ref{bifurcation equation z2y}) is expressed as 
\begin{equation}\label{bifurcation surface z2y}
    h_2{\alpha_1}^2+h_3{\alpha_2}^2+h_4{\alpha_3}^2+h_5e^2+\frac{1}{2\nu_0}-\frac{\nu_0}{2}=0.
\end{equation}
When $\alpha_3\in \mathbb{R}$, (\ref{bifurcation surface z2y}) describes a one-sheet hyperboloid in the $(\alpha_1,\alpha_2,\alpha_3)$ coordinates system, as shown in Fig.~\ref{fig6}(a). When $\alpha_3\in \sqrt{-1}\mathbb{R}$, replacing $\alpha_3$ with its imaginary part $\alpha_3/\sqrt{-1}$, (\ref{bifurcation surface z2y}) describes an ellipsoid, as illustrated in Fig.~\ref{fig6}(b).
Setting $\alpha_3=0$, the critical surface (\ref{bifurcation surface z2y}) reduces to a critical curve in the $(\alpha_1,\alpha_2)$ plane, as shown in Fig.~\ref{fig7}(a), which illustrates the orbital bifurcation in center manifolds associated with the breaking of the $S_2$-type symmetry. When the pair $(\alpha_1,\alpha_2)$ lies outside the critical curve, (\ref{bifurcation equation z2y}) has no non-zero solution. In this case, the expansion (\ref{formal sol}) describes non-bifurcated Lyapunov orbits and Lissajous orbits. On the other hand, when $(\alpha_1,\alpha_2)$ lies inside the critical ellipse, a pair of feasible solutions to the bifurcation equation exists, given by
\begin{equation}
    \eta=\pm\sqrt{- \frac{1}{h_1{\alpha_2}^2}(h_2{\alpha_1}^2+h_3{\alpha_2}^2+h_5e^2+\frac{1}{2\nu_0}-\frac{\nu_0}{2})},
\end{equation}
which describes two families of axial orbits and their corresponding quasi-axial orbits. Specifically, two families of axial orbits bifurcate from vertical Lyapunov periodic orbits under the critical condition
\begin{equation}
    |\alpha_2|\leq\sqrt{-\frac{1}{h_3}(h_5e^2+\frac{1}{2\nu_0}-\frac{\nu_0}{2})}.
\end{equation}
The bifurcated axial/quasi-axial orbits in the Sun-Earth system are illustrated in Fig.~\ref{fig7}(b).
\subsubsection{Case of coupling the motion in the $y$-direction to the motion in the $z$-direction.}
\begin{figure}[!htb]
    \centering
    \subfigure[]{
    \includegraphics[width=0.4\textwidth]{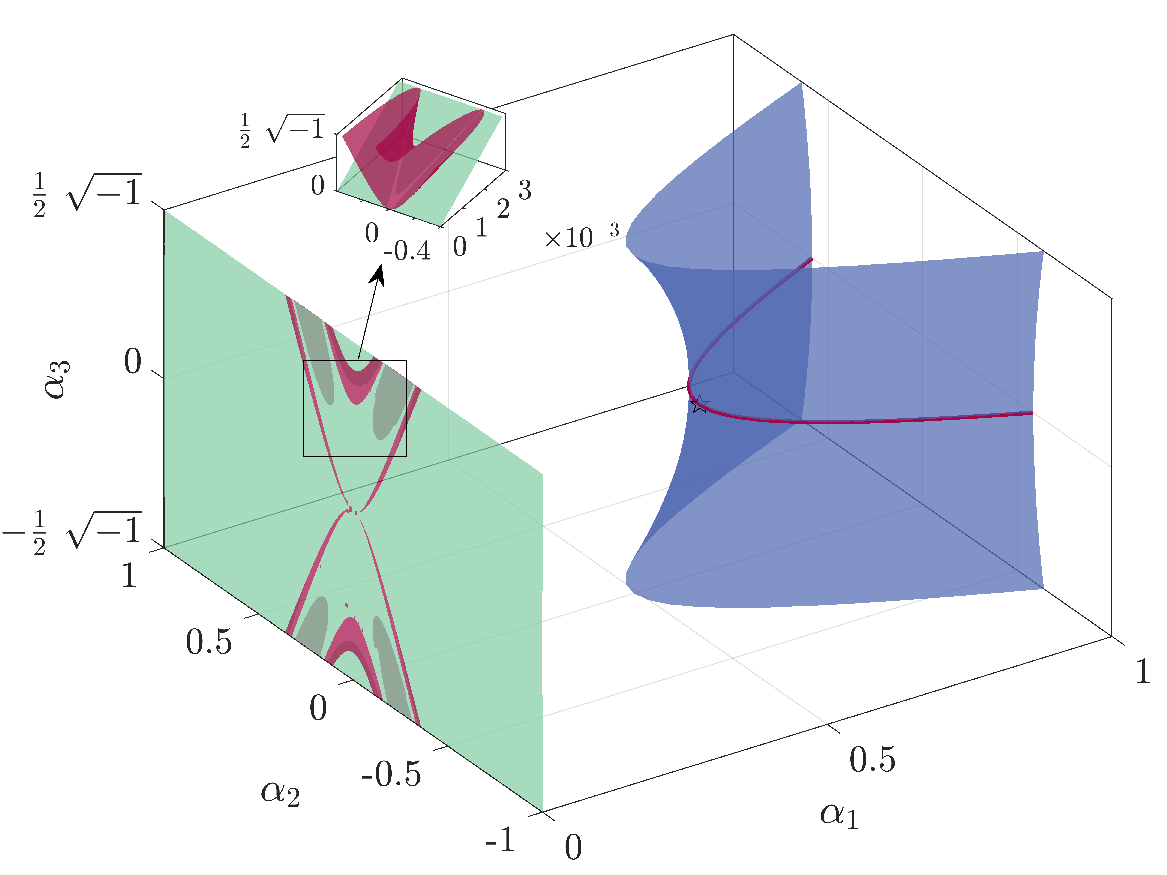}
    }
   \centering
    \subfigure[]{
    \includegraphics[width=0.4\textwidth]{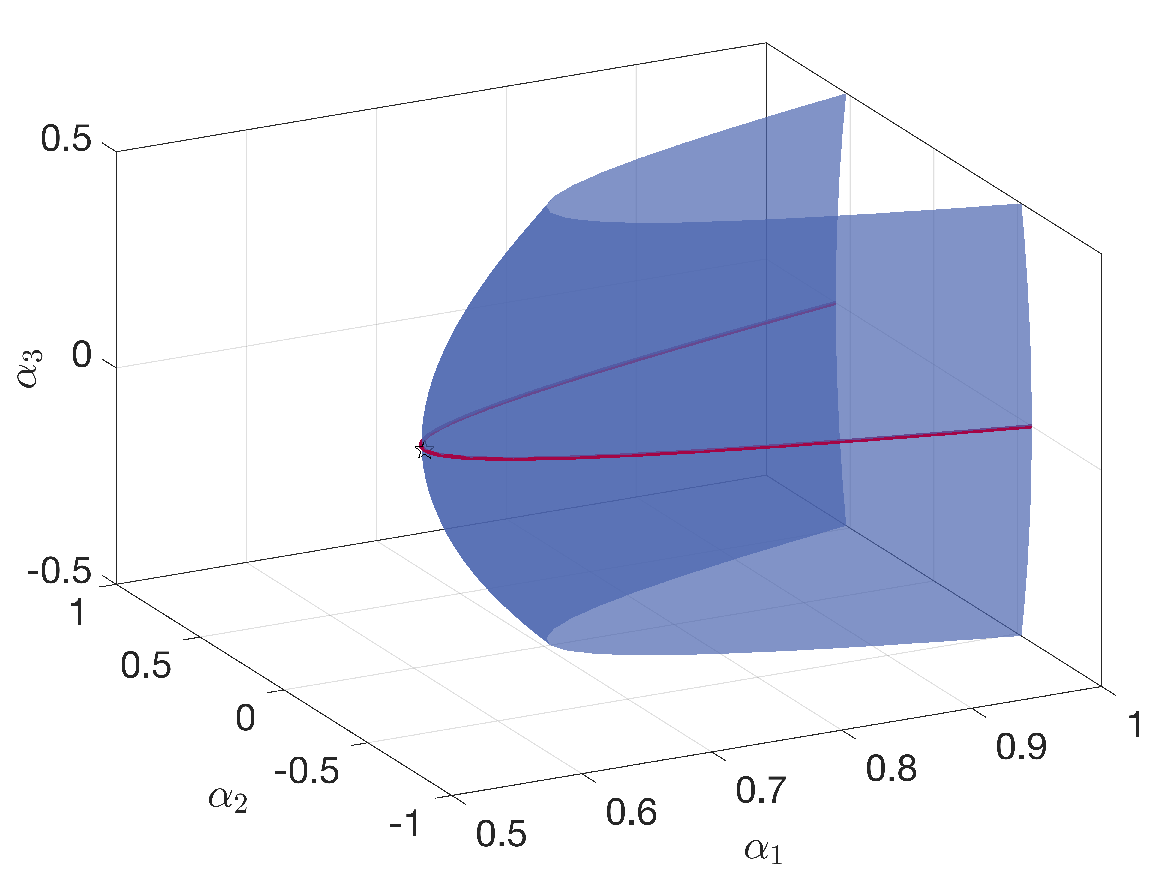}
    }
    \caption{Critical surfaces determined by the bifurcation equation for the case of considering a coupling effect from the motion in $y$-axis direction to the motion in the $z$-axis direction. (a) Critical surface with $\alpha_3$ associated with transit orbits. (b) Critical surface with $\alpha_3$ associated with non-transit orbits.}
    \label{fig8}
\end{figure}

It's obtained that axial/quasi-axial orbits can also bifurcate from planar Lyapunov periodic/Lissajous orbits. To demonstrate this, we consider the modification equations (\ref{modified ERTBP2}) which couples the motion in the $y$-direction to the motion in the $z$-direction. In this case, the bifurcation equation $\Delta=0$ is initially formulated by
\begin{equation}\label{y2z_bifurcation_eq}
    \Delta=(k_1{\alpha_1}^2+k_2{\alpha_3}^2)\eta^4+(k_3{\alpha_1}^2+k_4{\alpha_2}^2+k_5{\alpha_3}^2)\eta^2+(k_6{\alpha_1}^2+k_7{\alpha_2}^2+k_8{\alpha_3}^2+k_9e^2+\frac{{\nu_0}^2-{\omega_0}^2}{\kappa_1})=0 ,
\end{equation}
where $k_i(1\leq i\leq 9)$ depends solely on the system parameter $\mu$. The critical surfaces in this case are illustrated in Fig.~\ref{fig8}, similar to the case of coupling the motion in the $x$-direction to the $z$-direction.
\begin{figure}[!htb]
    \centering
    \subfigure[]{
    \includegraphics[width=0.4\textwidth]{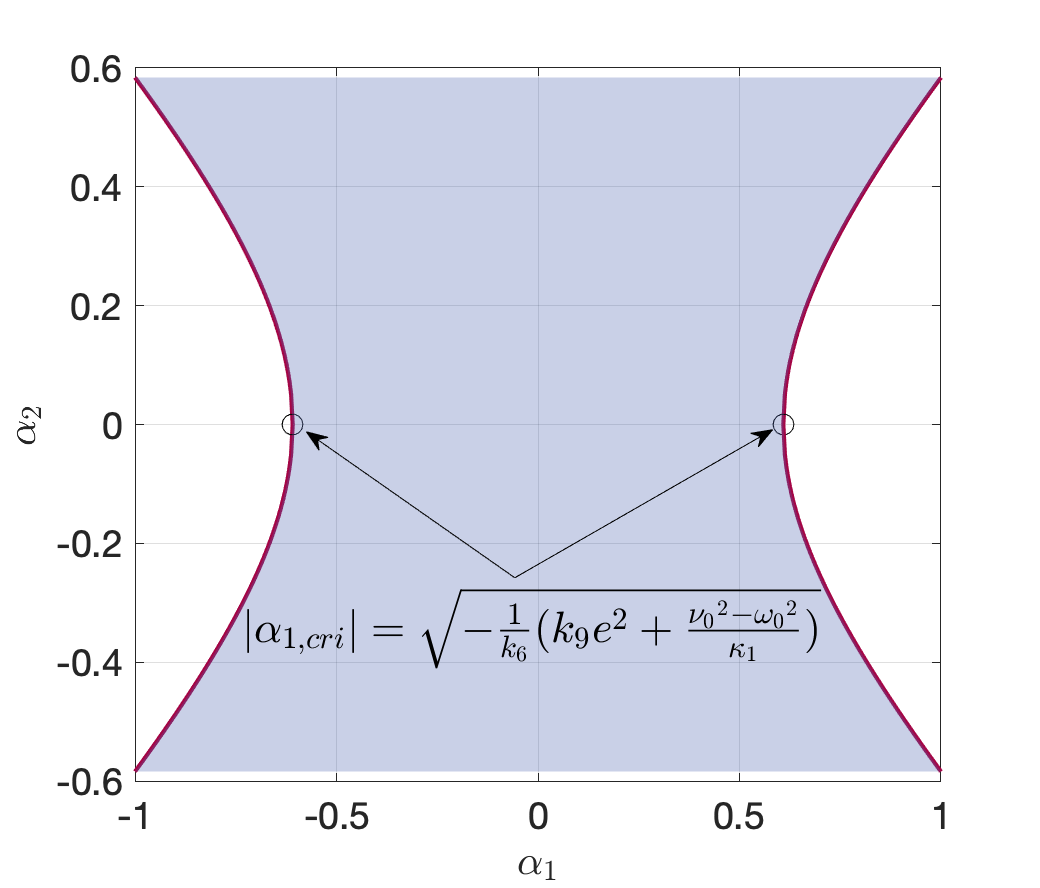}
    }
    \subfigure[]{
    \includegraphics[width=0.4\textwidth]{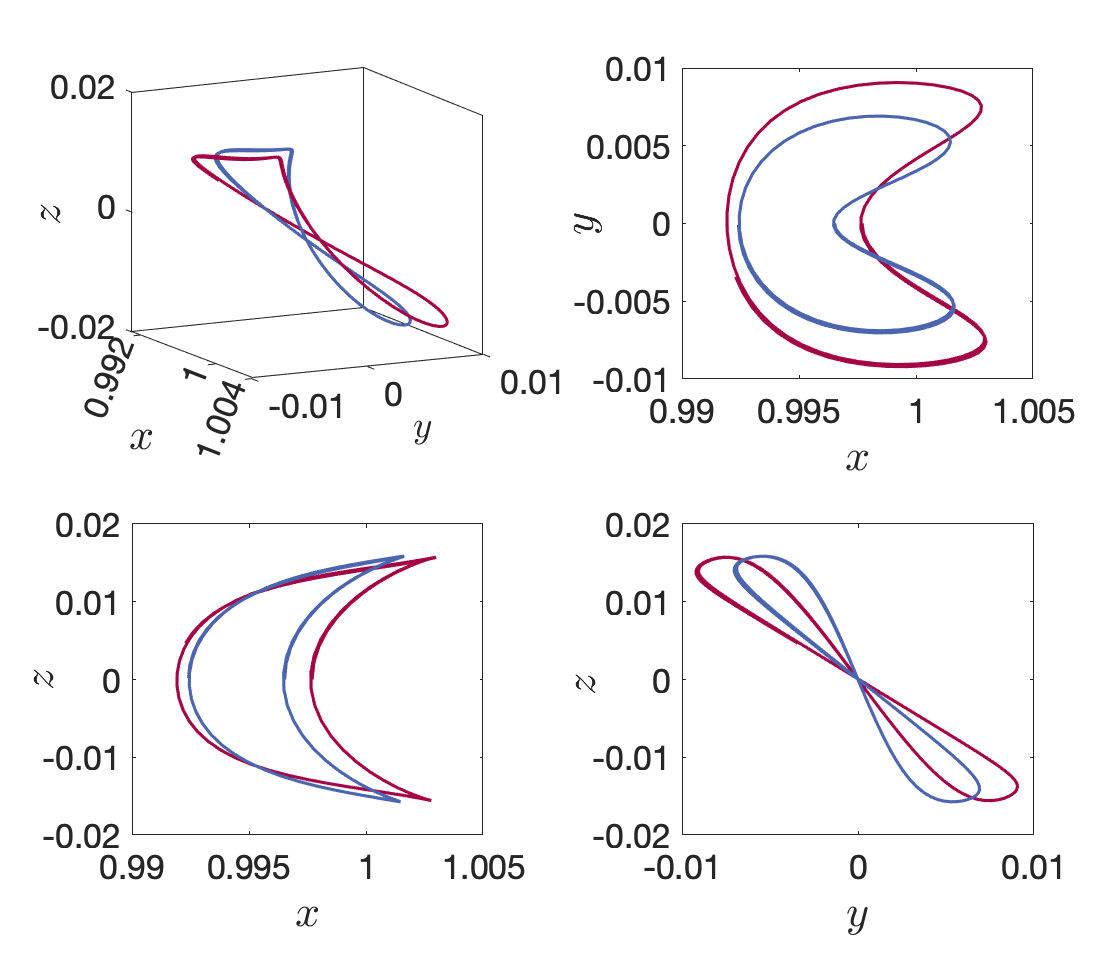}
    }
    \caption{(a): Bifurcation diagram of the third-order series solution of orbits in center manifolds in the case of coupling the motion in the $y$-direction to the motion in the $z$-direction. (b): Red: Axial orbit with $\alpha_1=0.28$. Blue: Axial orbit with $\alpha_1=0.2$.}
    \label{fig9}
\end{figure}

By letting the amplitude that associates with hyperbolic part of the solution be zero, the bifurcation equation (\ref{y2z_bifurcation_eq}) is reduced to
\begin{equation}\label{bifurcation equation center y2z}
    \Delta=k_1{\alpha_1}^2\eta^4+(k_3{\alpha_1}^2+k_4{\alpha_2}^2)\eta^2+(k_6{\alpha_1}^2+k_7{\alpha_2}^2+k_9e^2+\frac{{\nu_0}^2-{\omega_0}^2}{\kappa_1})=0.
\end{equation}
This equation defines a bifurcation curve in the $(\alpha_1,\alpha_2)$ plane, as illustrated in Fig.~\ref{fig9}(a). Similarly, feasible solutions $\eta$ to the bifurcation equation (\ref{bifurcation equation center y2z}) describe axial and quasi-axial orbits. In particular, when it satisfies that $\alpha_2=\alpha_3=0$, and the absolute value of $\alpha_1$ is smaller than some critical value given by 
\begin{equation}
    |\alpha_{1,cri}|=\sqrt{-\frac{1}{k_6}(k_9e^2+\frac{{\nu_0}^2-{\omega_0}^2}{\kappa_1})},
\end{equation}
there exists a pair of feasible solutions $\eta$, satisfying
\begin{equation}
    \eta=\pm\sqrt{\frac{-k_3{\alpha_1}^2-\sqrt{{k_3}^2{\alpha_1}^4-4k_1{\alpha_1}^2(k_6{\alpha_1}^2+k_9e^2+\frac{{\nu_0}^2-{\omega_0}^2}{\kappa_1})}}{2k_1{\alpha_1}^2}}.
\end{equation}
These bifurcated solutions describe two families of axial orbits that bifurcate from planar Lyapunov orbits as shown in Fig.~\ref{fig9}(b).

\section{Conclusions}\label{sec5}
In this paper, we present a semi-analytical framework to describe pitchfork bifurcations and symmetry breaking near collinear libration points in the ERTBP. By reformulating the hyperbolic components of the general solution within the complex field, we develop a unified trigonometric series-based approach to systematically exploit the inherent symmetries of the ERTBP. With the aid of coupling-induced bifurcation mechanisms, we achieve controlled symmetry breaking from the semi-analytical structures of non-bifurcated orbits. This is accomplished by constructing different bifurcation equations and introducing corresponding coupling coefficients that parameterize the transition between symmetric and asymmetric dynamical configurations.
A quantitative study of these parameterized bifurcation equations reveals comprehensive insights into the first-level bifurcation dynamics triggered from pitchfork bifurcations in the non-autonomous ERTBP.
Specifically, we demonstrate that the emergence of bifurcations is governed by the existence of feasible solutions $(\eta\not=0)$ to the bifurcation equation $\Delta=0$, where solutions are determined solely by the orbital eccentricity $e$ and three amplitude parameters $\alpha_i$ $(i=1,2,3)$.
These bifurcations encompass not only periodic/quasi-periodic orbits but also transit/non-transit orbits, thereby unifying the characterization of bifurcations of both central and hyperbolic dynamical behaviors.

Due to inherent limitations of local perturbation methods, semi-analytical frameworks within the Lindstedt-Poincaré method remain challenging for analyzing the “next-level” bifurcations, particularly for bifurcated orbits with large-amplitude excursions from libration points. Meanwhile, the trigonometric framework in representing the hyperbolic motion meets obstacles in providing a high-order solution of stable/unstable manifolds associated with bifurcated halo/axial orbits. Further exploration of methods for accomplishing systematic symmetry breaking from non-bifurcated structures is needed in future work.\\ 

\pdfbookmark[1]{Acknowledgements}{sec:acknowledgements}
\section*{Acknowledgements}
\noindent
The authors are grateful to Prof. Hayato Chiba for useful suggestions and comments. The author H.S. was supported by JST CREST, Grant Number JPMJCR2014, Japan.
The author M.L. was supported by the Sustainability Open Knowledge-Action Program by Connecting Multi-stakeholder Projects from Tohoku University and JST Moonshot, Grant Number JPMJMS2023, Japan.

\pdfbookmark[1]{Author Contributions}{sec:author_contributions}
\section*{Author Contributions}
\noindent
H.S. contributed to the methodology, software, validation,
formal analysis, investigation, and writing-original draft visualization. 
M.L. contributed to the conceptualization, software, writing-review and editing, supervision, project administration, and funding acquisition.

\pdfbookmark[1]{Data Availability}{sec:data_availability}
\section*{Data Availability}
\noindent
The datasets generated during and analyzed during the current study are available from the corresponding author upon reasonable request.

\pdfbookmark[1]{Declarations}{sec:declaration}
\section*{Declarations}
\noindent
Conflict of interest The authors declare no conflict of interest.

\pdfbookmark[1]{Appendix}{sec:appendix}
\section*{Appendix}
\label{appendix}
A detailed computation of the undetermined coefficients associated with breaking of the $S_1$-type symmetry is provided. According to the discussions in Section \ref{sec3}, initializing with the modified linear solution (\ref{modified linear 1}), high-order solutions of halo/quasi-halo orbits and their corresponding transit/non-transit orbits can be obtained through a regular perturbation procedure.
Starting with the solution up to order $n-1$ ($n\geq 2$), we substitute it into the modified model of (\ref{ERTBP3}) according to the modification (\ref{modified_ERTBP}) to obtain all known terms. Here, the first and second derivatives of $x$ with respect to $f$ are given by
$$
\begin{aligned}
    x'=\ &\omega\frac{\partial x}{\partial \theta_1}+\nu\frac{\partial x}{\partial \theta_2}+\sqrt{-1}\lambda\frac{\partial x}{\partial \theta_3}+\frac{\partial x}{\partial f},\\
    x''=\ &\omega^2\frac{\partial^2 x}{\partial {\theta_1}^2}+\nu^2\frac{\partial^2 x}{\partial {\theta_2}^2}-\lambda^2\frac{\partial^2 x}{\partial {\theta_3}^2}+\frac{\partial ^2x}{\partial f^2}+2\omega\nu\frac{\partial^2 x}{\partial \theta_1\partial \theta_2}+2\sqrt{-1}\omega\lambda\frac{\partial^2 x}{\partial \theta_1\partial \theta_3}+2\sqrt{-1}\nu\lambda\frac{\partial^2 x}{\partial \theta_2\partial \theta_3}\\
    &+2\omega\frac{\partial^2 x}{\partial\theta_1\partial f}+2\nu\frac{\partial^2 x}{\partial\theta_2\partial f}+2\sqrt{-1}\lambda\frac{\partial^2 x}{\partial\theta_3\partial f},
\end{aligned}
$$
while the derivatives of other coordinates can be calculated analogously.
These known terms refined from the left-hand side of the modified model are moved to the right-hand side and denoted by three new series $p,\ q\ \text{and}\ \tau$. 
Specifically, their $n$-th order terms $p_{ijkm}^{stur}, q_{ijkm}^{stur}$, and $\tau_{ijkm}^{stur}$ with $i+j+k+m=n$ are fully determined. In this way, the unknown coefficients of $n$-th order can be obtained by solving a sequence of linear equations, where the unknown terms on the left-hand side match the known terms on the right-hand side in (\ref{modified_ERTBP}).
Apart from three singular cases ($|i|+|j|+|k|=1$), the general form of these linear equations for solving undetermined coefficients in (\ref{formal sol}) is as follows:
\begin{equation}\label{eq1}
\begin{bmatrix}
    -{\tilde{\omega}_{stur}}^2-1-2c_2& -2\tilde{\omega}_{stur} \\
    -2\tilde{\omega}_{stur}& -{\tilde{\omega}_{stur}}^2+c_2-1 
\end{bmatrix}
\begin{bmatrix}
    x_{ijkm}^{stur}\\
    y_{ijkm}^{stur}
\end{bmatrix}
=
\begin{bmatrix}
    p_{ijkm}^{stur}\\
    q_{ijkm}^{stur}
\end{bmatrix},
\end{equation}
\begin{equation}\label{eq2}
\begin{aligned}
    (c_2-{\tilde{\omega}_{stur}}^2)z_{ijkm}^{stur}=\tau_{ijkm}^{stur}+\eta d_{0000}x_{ijkm}^{stur},
\end{aligned}
\end{equation}
where $\tilde{\omega}_{stur}=s\omega_0+t\nu_0+\sqrt{-1}u\lambda_0+r$ is a complex-valued constant. Due to the definitions of $\omega_0,\ \nu_0\ \text{and}\ \lambda_0$, the determinant of the coefficient matrix is non-zero. Hence, the undetermined coefficients in (\ref{eq1}) and (\ref{eq2}) can be solved directly. Nevertheless, when specific quartets $(s,t,u,r)$ are selected, singularities emerge in the linear equation system. These exceptional cases include the following scenarios:\\
\textbf{Case 1:} (\textit{s,t,u,r})$=(1,0,0,0).$\\
In this case, the linear algebraic equations are formulated as
\begin{equation}\label{sin_eq1}
\begin{bmatrix}
    -{\omega_0}^2-1-2c_2& -2\omega_0 \\
    -2\omega_0& -{\omega_0}^2+c_2-1 
\end{bmatrix}
\begin{bmatrix}
    x_{ijkm}^{1000}\\
    y_{ijkm}^{1000}
\end{bmatrix}
+\begin{bmatrix}
    -2(\kappa_1+\omega_0)\omega_{i-1jkm}\\
    -2(1+\kappa_1\omega_0)\omega_{i-1jkm}
\end{bmatrix}
=
\begin{bmatrix}
    p_{ijkm}^{1000}\\
    q_{ijkm}^{1000}
\end{bmatrix}
+
\begin{bmatrix}
    \Omega_{i-1jkm}\\
    \kappa_1\Omega_{i-1jkm}
\end{bmatrix},
\end{equation}
\begin{equation}\label{x2z_d_ijkm}
\begin{aligned}
    (c_2-{\omega_0}^2)z_{ijkm}^{1000}-2\eta\omega_0\omega_{i-1jkm}=\tau_{ijkm}^{1000}+\eta d_{0000}x_{ijkm}^{1000}+\eta d_{i-1jkm}+\eta\Omega_{i-1jkm},
\end{aligned}
\end{equation}
where $\Omega_{i-1jkm}=\omega^2-2\omega_{0}\omega_{i-1jkm}$ represents the corresponding known term of $\omega^2$ of order $n-1$. 
Equation (\ref{sin_eq1}) is singular. We set $x_{ijkm}^{1000}=0$, and then, $y_{ijkm}^{1000}$ and $\omega_{i-1jkm}$ can be obtained by dealing with the following reduced regular system of equations:
\begin{equation}
    \begin{bmatrix}
        -2\omega_0 & -2(\kappa_1+\omega_0)\\
        -{\omega_0}^2+c_2-1 & -2(1+\kappa_1\omega_0)
    \end{bmatrix}
    \begin{bmatrix}
        y_{ijkm}^{1000}\\
        \omega_{i-1jkm}
    \end{bmatrix}
    =
    \begin{bmatrix}
    p_{ijkm}^{1000}\\
    q_{ijkm}^{1000}
\end{bmatrix}
+
\begin{bmatrix}
    \Omega_{i-1jkm}\\
    \kappa_1\Omega_{i-1jkm}
\end{bmatrix},
\end{equation}
Similarly, by letting $z_{ijkm}^{1000}=0$, the coefficients $d_{i-1jkm}$ are solvable as well.\\
\textbf{Case 2:} (\textit{s,t,u,r})=$(0,1,0,0).$\\
The linear equations of the undetermined coefficients are as follows:
\begin{equation}
\begin{bmatrix}
    -{\nu_0}^2-1-2c_2& -2\nu_0 \\
    -2\nu_0& -{\nu_0}^2+c_2-1 
\end{bmatrix}
\begin{bmatrix}
    x_{ijkm}^{0100}\\
    y_{ijkm}^{0100}
\end{bmatrix}
=
\begin{bmatrix}
    p_{ijkm}^{0100}\\
    q_{ijkm}^{0100}
\end{bmatrix},
\end{equation}
\begin{equation}
\begin{aligned}
    (c_2-{\nu_0}^2)z_{ijkm}^{0100}-2\nu_0\nu_{ij-1km}=\tau_{ijkm}^{0100}+\eta d_{0000}x_{ijkm}^{0100}+\Lambda_{ij-1km}.
\end{aligned}
\end{equation}
In this case, $x_{ijkm}^{0100}$ and $y_{ijkm}^{0100}$ can be solved directly. $\Lambda_{ij-1km}$ represents the known term of $\nu^2$ of order $n-1$. Letting $z_{ijkm}^{0100}=0$, we have $\nu_{ij-1km}=-\left(\tau_{ijkm}^{0100}+\eta d_{0000}x_{ijkm}^{0100}+\Lambda_{ij-1km}\right)/{2\nu_0}$.\\
\textbf{Case 3:} (\textit{s,t,u,r})=$(0,0,1,0).$\\
In this case, the undetermined coefficients should satisfy
\begin{equation}
\begin{bmatrix}
    {\lambda_0}^2-1-2c_2& -2\sqrt{-1}\lambda_0 \\
    -2\sqrt{-1}\lambda_0& {\lambda_0}^2+c_2-1 
\end{bmatrix}
\begin{bmatrix}
    x_{ijkm}^{0010}\\
    y_{ijkm}^{0010}
\end{bmatrix}
+\begin{bmatrix}
    2(\lambda_0+\kappa_2)\lambda_{ijk-1m}\\
    2(\sqrt{-1}\kappa_2\lambda_0-\sqrt{-1})\lambda_{ijk-1m}
\end{bmatrix}
=
\begin{bmatrix}
    p_{ijkm}^{0010}\\
    q_{ijkm}^{0010}
\end{bmatrix}
+
\begin{bmatrix}
    \Gamma_{ijk-1m}\\
    \sqrt{-1}\kappa_2\Gamma_{ijk-1m}
\end{bmatrix},
\end{equation}
\begin{equation}
\begin{aligned}
    (c_2+{\lambda_0}^2)z_{ijkm}^{0010}+2\eta\kappa_3\lambda_0\lambda_{ijk-1m}=\tau_{ijkm}^{0010}+\eta d_{0000}x_{ijkm}^{0010}+\eta d_{ijk-1m}+\eta\kappa_3\Gamma_{ijk-1m},
\end{aligned}
\end{equation}
where $\Gamma_{ijk-1m}$ represents the known term of $\lambda^2$ of order $n-1$.
Like the \textbf{Case 1}, by setting $x_{ijkm}^{0010}$ to zero, the remaining terms $y_{ijkm}^{0010},\ z_{ijkm}^{0010},\ \text{and}\ \lambda_{ijk-1m}$ can be computed correspondingly. More precisely, the reduced equations can be represented as
\begin{equation}
    \begin{bmatrix}
        -2\sqrt{-1}\lambda_0 & 2(\lambda_0+\kappa_2)\\
        {\lambda_0}^2+c_2-1 & 2(\sqrt{-1}\kappa_2\lambda_0-\sqrt{-1})
    \end{bmatrix}
    \begin{bmatrix}
        y_{ijkm}^{0010}\\
        \lambda_{ijk-1m}
    \end{bmatrix}
    =
\begin{bmatrix}
    p_{ijkm}^{0010}\\
    q_{ijkm}^{0010}
\end{bmatrix}
+
\begin{bmatrix}
    \Gamma_{ijk-1m}\\
    \sqrt{-1}\kappa_2\Gamma_{ijk-1m}
\end{bmatrix},
\end{equation}
\begin{equation}
    \begin{aligned}
        (c_2+{\lambda_0}^2)z_{ijkm}^{0010}=-2\eta\kappa_3\lambda_0\lambda_{ijk-1m}+\tau_{ijkm}^{0010}+\eta d_{ijk-1m}+\eta\kappa_3\Gamma_{ijk-1m}.
    \end{aligned}
\end{equation}
The corresponding undetermined coefficients can be obtained by solving the equations above.

\addcontentsline{toc}{section}{References}
\bibliography{sample}

\begin{thebibliography}{26}
\newcommand{\enquote}[1]{``#1''}
\providecommand{\natexlab}[1]{#1}
\providecommand{\url}[1]{\texttt{#1}}
\providecommand{\urlprefix}{URL }
\expandafter\ifx\csname urlstyle\endcsname\relax
  \providecommand{\doi}[1]{\discretionary{}{}{}https://doi.org/#1}\else
  \providecommand{\doi}[1]{\discretionary{}{}{}\urlstyle{rm}\url{https://doi.org/#1}}\fi

\bibitem[{Gómez and Mondelo(2001)}]{GomezMondelo}
Gómez, G., and Mondelo, J.~M., \enquote{The dynamics around the collinear equilibrium points of the RTBP,} \emph{Physica D: Nonlinear Phenomena}, Vol. 157, No.~4, 2001, pp. 283--321.
\newblock \doi{10.1016/S0167-2789(01)00312-8}.

\bibitem[{Gómez et~al.(2003)Gómez, Masdemont, and Mondelo}]{GomezMasdemontMondelo}
Gómez, G., Masdemont, J.~J., and Mondelo, J.~M., \emph{Libration point orbits: a survey from the dynamical point of view}, World-Scientific, Singapore, 2003, pp. 311--372.
\newblock \doi{10.1142/9789812704849_0016}.

\bibitem[{Broucke(1969)}]{Broucke}
Broucke, R., \enquote{Stability of periodic orbits in the elliptic, restricted three-body problem,} \emph{AIAA Journal}, Vol.~7, No.~6, 1969, pp. 1003--1009.
\newblock \doi{10.2514/3.5267}.

\bibitem[{Ovenden and Roy(1961)}]{Ovenden}
Ovenden, M.~W., and Roy, A.~E., \enquote{On the use of the Jacobi integral of the restricted three-body problem,} \emph{Monthly Notices of the Royal Astronomical Society}, Vol. 123, 1961, pp. 1--14.
\newblock \doi{10.1093/mnras/123.1.1}.

\bibitem[{Parker and Anderson(2014)}]{parker2014low}
Parker, J.~S., and Anderson, R.~L., \emph{Low-energy lunar trajectory design}, John Wiley \& Sons, Ltd, Hoboken, 2014.
\newblock \doi{10.1007/s10569-012-9457-4}.

\bibitem[{Peng et~al.(2017)Peng, Bai, Masdemont, Gómez, and Xu}]{Peng_2}
Peng, H., Bai, X.~L., Masdemont, J.~J., Gómez, G., and Xu, S.~J., \enquote{Libration transfer design using patched elliptic three-body models and graphics processing units,} \emph{Journal of Guidance, Control, and Dynamics}, Vol.~40, No.~12, 2017, pp. 3155--3166.
\newblock \doi{10.2514/1.G002692}.

\bibitem[{Jorba-Cuscó and Epenoy(2021)}]{Jorba_Cusco}
Jorba-Cuscó, M., and Epenoy, R., \enquote{Low-fuel transfers from Mars to quasi-satellite orbits around Phobos exploiting manifolds of tori,} \emph{Celest Mech Dyn Astr}, Vol. 133, No.~20, 2021.
\newblock \doi{10.1007/s10569-021-10017-9}.

\bibitem[{Shirobokov et~al.(2017)Shirobokov, Trofimov, and Ovchinnikov}]{Shirobokov}
Shirobokov, M., Trofimov, S., and Ovchinnikov, M., \enquote{Survey of station-keeping techniques for libration point orbits,} \emph{Journal of Guidance, Control, and Dynamics}, Vol.~40, No.~5, 2017, pp. 1085--1105.
\newblock \doi{10.2514/1.G001850}.

\bibitem[{Gurfil and Meltzer(2007)}]{Gurfil}
Gurfil, P., and Meltzer, D., \enquote{Semi-analytical method for calculating the elliptic restricted three-body problem monodromy matrix,} \emph{Journal of Guidance, Control, and Dynamics}, Vol.~30, No.~1, 2007, pp. 266--271.
\newblock \doi{10.2514/1.22871}.

\bibitem[{Antoniadou and Voyatzis(2013)}]{Antoniadou}
Antoniadou, K.~I., and Voyatzis, G., \enquote{2/1 resonant periodic orbits in three dimensional planetary systems,} \emph{Celest Mech Dyn Astr}, Vol. 115, 2013, pp. 161--184.
\newblock \doi{10.1007/s10569-012-9457-4}.

\bibitem[{Peng and Xu(2015)}]{Peng_1}
Peng, H., and Xu, S., \enquote{Stability of two groups of multi-revolution elliptic halo orbits in the elliptic restricted three-body problem,} \emph{Celest Mech Dyn Astr}, Vol. 123, 2015, pp. 279--303.
\newblock \doi{10.1007/s10569-015-9635-2}.

\bibitem[{Ferrari and Lavagna(2018)}]{Ferrari}
Ferrari, F., and Lavagna, M., \enquote{Periodic motion around libration points in the elliptic restricted three-body problem,} \emph{Nonlinear Dyn}, Vol.~93, No.~6, 2018, p. 453–462.
\newblock \doi{10.1007/s11071-018-4203-4}.

\bibitem[{Paez and Guzzo(2021)}]{Paez_Guzzo}
Paez, R.~I., and Guzzo, M., \enquote{Transits close to the Lagrangian solutions $L_1$, $L_2$ in the elliptic restricted three-body problem,} \emph{Nonlinearity}, Vol.~34, No.~9, 2021, pp. 6417--6449.
\newblock \doi{10.1088/1361-6544/ac13be}.

\bibitem[{Jorba et~al.(2024)Jorba, Nicolás, and Rodríguez}]{Jorba_chaos}
Jorba, A., Nicolás, B., and Rodríguez, O., \enquote{A dynamical study of Hilda asteroids in the circular and elliptic RTBP,} \emph{Chaos}, Vol.~34, No.~12, 2024.
\newblock \doi{10.1063/5.0234410}.

\bibitem[{Farquhar(1969)}]{Farquhar}
Farquhar, R.~W., \enquote{The control and use of libration-point satellites,} Ph.D. thesis, Stanford Univ., Stanford, CA, 1969.

\bibitem[{Richardson(1980)}]{Richardson}
Richardson, D.~L., \enquote{Analytic construction of periodic orbits about the collinear points,} \emph{Celestial mechanics}, Vol.~20, No.~3, 1980, pp. 241--253.
\newblock \doi{10.1007/BF01229511}.

\bibitem[{Masdemont(2005)}]{Masdemont}
Masdemont, J.~J., \enquote{High-order expansions of invariant manifolds of libration point orbits with applications to mission design,} \emph{Dynamical Systems}, Vol.~20, No.~1, 2005, pp. 59--113.
\newblock \doi{10.1080/14689360412331304291}.

\bibitem[{Lei et~al.(2013)Lei, Xu, Hou, and Sun}]{Lei}
Lei, H.~L., Xu, B., Hou, X.~Y., and Sun, Y.~S., \enquote{High-order solutions of invariant manifolds associated with libration point orbits in the elliptic restricted three-body system,} \emph{Celest Mech Dyn Astr}, Vol. 117, No.~4, 2013, pp. 349--384.
\newblock \doi{10.1007/s10569-013-9515-6}.

\bibitem[{Jorba and Masdemont(1998)}]{Jarba_Masdemont}
Jorba, A., and Masdemont, J.~J., \enquote{Dynamics in the centre manifold of the collinear points of the restricted three body problem,} \emph{Physica D: Nonlinear Phenomena}, Vol. 132, 1998.
\newblock \doi{10.1016/S0167-2789(99)00042-1}.

\bibitem[{Paez and Guzzo(2022)}]{Paez_Guzzo_Ana}
Paez, R.~I., and Guzzo, M., \enquote{On the semi-analytical construction of halo orbits and halo tubes in the elliptic restricted three-body problem,} \emph{Physica D: Nonlinear Phenomena}, Vol. 439, 2022.
\newblock \doi{10.1016/j.physd.2022.133402}.

\bibitem[{Celletti et~al.(2024)Celletti, Lhotka, and Pucacco}]{Celleti}
Celletti, A., Lhotka, C., and Pucacco, G., \enquote{The dynamics around the collinear points of the elliptic three-body problem: A normal form approach,} \emph{Physica D: Nonlinear Phenomena}, Vol. 468, 2024.
\newblock \doi{10.1016/j.physd.2024.134302}.

\bibitem[{Lin and Chiba(2025)}]{Lin_1}
Lin, M., and Chiba, H., \enquote{Bifurcation mechanism of quasi-Halo orbit from Lissajous orbit,} \emph{Journal of Guidance, Control, and Dynamics}, Vol.~48, No.~1, 2025.
\newblock \doi{10.2514/1.G008233}.

\bibitem[{Lin et~al.(2024)Lin, Luo, and Chiba}]{Lin_2}
Lin, M., Luo, T., and Chiba, H., \enquote{Semi-analytical computation of bifurcation of orbits near collinear libration point in the restricted three-body problem,} \emph{Physica D: Nonlinear Phenomena}, Vol. 470, 2024.
\newblock \doi{10.1016/j.physd.2024.134404}.

\bibitem[{Doedel et~al.(2007)Doedel, Nomanov, Paffenroth, Keller, Dichmann, Galán-Vioque, and Vanderbauwhede}]{Doedel}
Doedel, E.~J., Nomanov, V.~A., Paffenroth, R.~C., Keller, H.~B., Dichmann, D.~J., Galán-Vioque, J., and Vanderbauwhede, A., \enquote{Elemental periodic orbits associated with the libration points in the circular restricted 3-body problem,} \emph{International Journal of Bifurcation and Chaos}, Vol.~17, No.~08, 2007, pp. 2625--2677.
\newblock \doi{10.1142/S0218127407018671}.

\bibitem[{Szebehely(1967)}]{Szebehely}
Szebehely, V., \enquote{Chapter 10 - Modifications of the restricted problem,} \emph{Theory of orbit}, Academic Press, 1967, pp. 556--652.
\newblock \doi{10.1016/B978-0-12-395732-0.50016-7}.

\bibitem[{Celletti(2010)}]{Celletti_book}
Celletti, A., \emph{Stability and chaos in celestial mechanics}, Springer Berlin, Heidelberg, 2010.
\newblock \doi{10.1007/978-3-540-85146-2}.

\end{thebibliography}

\end{document}